\documentclass[onecolumn]{autart}    
\usepackage{indentfirst}
\usepackage{amsmath,amssymb,amsfonts}
\usepackage{algorithmic}
\usepackage{graphicx}
\usepackage{textcomp}
\usepackage{epstopdf}
\usepackage{bm}

\usepackage{soul}
\usepackage{graphicx}          
\usepackage[dvips]{epsfig}    

\begin{document}

\begin{frontmatter}

\title{On Parameter Selection of Nonsingular Predefined-time Terminal Sliding Mode With the Fixed-time Convergence Guarantee} 

\thanks[footnoteinfo]{Corresponding author Wen Yan. Tel. 18683977138.}

\author[Paestum]{Wen Yan}\ead{essay.yan@qq.com},    

\address[Paestum]{College of Electrical Engineering, Sichuan University, Chengdu 610065, China}  


\begin{abstract}                          
This paper study the parameter selection of predefined-time sliding mode and try to design a general nonsingular predefined-time terminal sliding mode. 1). On parameter selection: Some existing predefined-time sliding modes are designed to focus on the reaching time and ignore the characterization of the equilibrium point i.e. $\bm x_e= \bm 0$. In actual engineer, the system can only converge to a small neighborhood near the equilibrium point because of the existence of uncertainty i.e. $\bm x_e \to \bm 0$. The actual equilibrium point should be deduced by taking the limit but not directly solve the equilibrium point. Hence, the actual selection of exponential term is suggested to be not 1 to make system convergence to the equilibrium point within predefined time. 2). On singularity-avoidance: Based on the exponential feature of predefined-time stability systems, a mathematic concept of switching sliding mode \cite{su1994adaptive} is applied in the design of nonsingular predefined-time terminal sliding mode. However, this class of sliding mode switching method will cause a singular problem, which can be described as control input $\bm u \to \infty$ if the system state $\bm x \to \bm 0$ and any-order time derivative of $\bm x$ is not 0. Hence, a non-singular selection condition for each recursive sliding mode surface is explored in this paper to realize the global finiteness of control input. In addition, an unnecessary chattering of control input will also be caused by the sliding mode switching frequency, and a novel switching condition will be applied in the controller to reduce this chattering. Simulations will be carried out to validate the effect of the novel nonsingular predefined-time sliding mode control method to some extent.
\end{abstract}

\end{frontmatter}



\section{Introduction}

Sliding mode control is widely studied by researchers for its robustness, time convergence and high-precision \cite{young1999control,benchaib1999real,yu2005continuous}. An important research direction of sliding mode control is finite time control \cite{ding2011survey}. With the development of finite time sliding mode control, the characteristic of a terminal attractor is explored gradually to obtain faster convergence rate of system \cite{venkataraman1993control}. It is not satisfied that the convergence time depends on the initial state, a fixed-time attractivity is presented in \cite{polyakov2011nonlinear} to obtain a fixed-time stability system. However, the upper bound of the settling time of traditional fixed-time is usually determined by two or more parameters, hence, a simplifying fixed-time stability system is presented in \cite{sanchez2015predefined} to make the upper bound of settling time explicit in the system function. In addition, the improved fixed-time stability system can also obtain an upper bound of settling time without enlarging and reducing mathematically \cite{jimenez2017semi,jimenez2018semi}, but most of traditional fixed-time stability systems enlarged the upper bound of settling time. Especially, a note on predefined-time stability in \cite{jimenez2018note} pointed that the predefined-time stability of system can be guaranteed without the exponential term. This result allows a construction of the predefined-time stabilizing controller with only polynominal term (instead of the exponential one) and its Lyapunov characterization analysis is perfect, but this class of predefined-time stability function \cite{jimenez2018note,sanchez2016non,aldana2019design,jimenez2017predefined} may not consider if system can reach and stay in the neighborhood of equilibrium point within predefined time. The related researchers study a reconstruction method of controller in Definition 3.1 of \cite{sanchez2015predefined} to make system convergence to the equilibrium point theoretically, but it seems that the actual characterization of the equilibrium point will be different from the description in \cite{sanchez2015predefined} by considering the analysis and simulation in this paper.

Furthermore, it is worth noting that almost all fixed-time or predefined-time sliding mode controllers have a singular problem. This singularity can be defined as a case where the control input may be infinite if the system state $\bm x \to \bm 0$ and any-order time derivative of $\bm x$ is not 0 or $\infty$. For a second-order system, related researchers have presented many non-singulary control methods, and these methods can be divided into two categories: one is to transfer the singularity caused by the state $\bm x=0$ to the time derivative of the state $\bm \dot x=0$ \cite{zuo2014non,sanchez2016non,huang2017fixed}, another one is to deal with singular term of controller directly \cite{ni2016fast,corradini2018nonsingular,ni2019predefined}. The first category seems to deal well with the singularity proposed in this paper, but the control input may be infinite if $\bm \dot x=0$. The value of this time derivative is possible in practice. Another one(except the method in \cite{ni2019predefined}) can only guarantee finite-time stability of system, because its essence is to limit the max value of control input while ignoring the infinite control input required by the controller to ensure the fixed time convergence of the system when system is in singularity, then the actual convergence time may be longer than the predefined maximum convergence time. Fortunately, based on the exponential feature of predefined-time stability systems, a traditional sliding mode switching strategy \cite{su1994adaptive} can be applied to design a nonsingular predefined-time terminal sliding mode controller in \cite{ni2019predefined}. This controller can guarantee not only the strict fixed time stability of system but also the finiteness of the control input, however, it is difficult to find an appropriate switching sliding mode. For a higher-order system, there seems to be little research on non-singular fixed-time recursive sliding mode control. The above-mentioned methods may be too complex if apply them directly to higher-order system, hence, it is necessary to study a simplifying sidling mode switching scheme to overcome the singularity of fixed-time terminal sidling mode control method for second-order or high-order system.

The contribution of this paper lies on: 1). The relationship between the non-time parameter value of predefined-time stability system and the characterization of equilibrium point is studied to explore the better value range of parameters. 2). A simplifying sidling mode switching scheme is studied to overcome the singularity of fixed-time terminal sliding mode control. 3). A nonsingular predefined-time sliding mode controller for higher-order system is studied to explore the max range of non-time parameter value.

\section{Preliminaries}

\textbf{Definition 1:} (see \cite{sanchez2015predefined}). Consider a system:
\begin{equation} \dot{\bm  x}(t) =  - \frac{1}{{{T_c}}}f(\bm x(t);m),~{\rm{ }}\bm x(0) = {\bm x_0}.\label{eq:1}\end{equation}
where $\bm x = {[{x_1},{x_2} \cdot  \cdot  \cdot {x_n}]^T} \in {R^n}$ is the vector of system state, $f(\bm x):D \to {R^n}$ is a continuous function and $f(\bm x) = 0$. The system \eqref{eq:1} is predefined-time stable if system is globally fixed-time stable and the minimum upper bound of settling time $T(\bm x_0)$ can be deduced without enlarging and reducing in math.

\textbf{Lemma 1:} (see \cite{jimenez2018note}). Consider the following system:
\begin{equation}  \dot {\bm x} =  - \frac{1}{{{T_c}}}{\left\| {\frac{{\partial W(\bm x)}}{{\partial \bm  x}}} \right\|^{ - 2}}{\left[ {\frac{{\partial W(\bm x)}}{{\partial \bm x}}} \right]^T}.\label{eq:2}\end{equation}
where $T_c>0$, and $W(\bm x)$ satisfies the following conditions:
\begin{equation}\begin{array}{l}
(i)0 < W(\bm x) < 1{\rm{ }}~for~{\rm{ }}~all~{\rm{ }}\bm x \in {R^n}\\
(ii)W(\bm x) = 0{\rm{ ~if ~and ~only~ if~ }}\bm x = 0\\
(iii)W(\bm x) \to 1{\rm{~ as~ }}\left\|\bm  x \right\| \to \infty \\
(iv)W(\bm x){\rm{ }}~is~{\rm{ }}{C^1}{\rm{~ in~ }}{R^n}\backslash \{ 0\} \\
(v)\left\| {\frac{{\partial W(\bm x)}}{{\partial \bm x}}} \right\| \ne 0{\rm{~ for~ all ~}}\bm x \in {R^n}
\end{array}\label{eq:3}\end{equation}
Then, the system can be said predefined-time stable, and the minimum upper bound of settling time is $T_c$.

\textbf{Lemma 2:} Consider the following system:
\begin{equation} \dot {\bm x}= - \frac{1}{{{T_c}}}{\left\| {\frac{{\partial W({\left\|\bm x \right\|^m})}}{{\partial \bm x}}} \right\|^{ - 2}}{\left[ {\frac{{\partial W({\left\|\bm x \right\|^m})}}{{\partial \bm x}}} \right]^T}.\label{eq:4}\end{equation}
where $W({\left\|\bm x \right\|^m})$ satisfies the conditions in \eqref{eq:3}. The exponential term should be emphasized in here and $0<m<1$ $(m \ne 1)$. The \eqref{eq:3}-(v) should be modified as $ {\frac{{\partial W({{\left\| \bm x \right\|}^m})}}{{\partial {{\bm x }}}}}  > 0{\rm{~ for ~all~}}\bm x \in {R^n}$. \textbf{In addition, a supplementary condition should be: $ {\frac{{\partial W({{\left\| \bm x \right\|}^m})}}{{\partial {{{{\left\| \bm x \right\|}^m} }}}}} \ne 0 {\rm{~for~all~}} \bm x \in {R^n}$.}

Proof: The system \eqref{eq:4} can be written as:
\begin{equation}\begin{array}{l}
\dot {\bm x} =  - \frac{1}{{{T_c}}}{\left\| {\frac{{\partial W({{\left\| \bm x \right\|}^m})}}{{\partial \bm x}}} \right\|^{ - 2}}{\left[ {\frac{{\partial W({{\left\|\bm x \right\|}^m})}}{{\partial \bm x}}} \right]^T}\\
~~ = - \frac{1}{{{T_c}}}\frac{{{{\left[ {\frac{{\partial W({{\left\|\bm x \right\|}^m})}}{{\partial \bm x}}} \right]}^T}\frac{{\partial W({{\left\|\bm x \right\|}^m})}}{{\partial\bm x}}}}{{{{\left\| {\frac{{\partial W({\bm x^m})}}{{\partial\bm x}}} \right\|}^{\rm{2}}}\frac{{\partial W({\bm x^m})}}{{\partial\bm x}}}}~~~~~~\\
~~= - \frac{1}{{{T_c}}}{\left[ {\frac{{\partial W({{\left\|\bm x \right\|}^m})}}{{\partial\bm x}}} \right]^{ - 1}}~~~~~~~~~~~~~~
\end{array}.\label{eq:5}\end{equation}
1)	By considering $m=1$, the $W({\left\| \bm x \right\|^m})$ still satisfies the condition in \eqref{eq:3}. However, a limit will be exampled in the following to emphasize $m \ne 1$ :
\begin{equation} \begin{array}{l}
\mathop {\lim }\limits_{\bm x \to 0} \dot {\bm x} =  - \frac{1}{{{T_c}}}{\left[ {\frac{{\partial W(\left\|\bm  x \right\|)}}{{\partial \left\|\bm  x \right\|}}\frac{{{\bm x^T}}}{{\left\| \bm x \right\|}}} \right]^{ - 1}}~~\\
~~~~~~{\rm{      }}\mathop  = \limits^{y = \left\| \bm x \right\|}  - \frac{1}{{{T_c}}}{\left[ {\frac{{\partial W(y)}}{{\partial y}}} \right]^{ - 1}}\frac{\bm x}{{\left\|\bm  x \right\|}}\\
{\rm{        }} = \left\{ {\begin{array}{*{20}{l}}
{ - \frac{1}{{{T_c}}}{{\left[ {\frac{{\partial W(y)}}{{\partial y}}} \right]}^{ - 1}}{\bm E_{n1}}{\rm{ ~~~~~~~~~~~~~~~~~~~~~~~~~~~~~~~~~~~~,}}\bm x \to {0^ + }}\\
{ \frac{1}{{{T_c}}}{\left[ {\frac{{\partial W(y)}}{{\partial y}}} \right]^{ - 1}}\left[ {\begin{array}{*{20}{c}}
{sign({x_1})}&0&0\\
0& \ddots &0\\
0&0&{sign({x_n})}
\end{array}} \right]{\bm E_{n1}}{\rm{,}}others}\\
~~{\frac{1}{{{T_c}}}{{\left[ {\frac{{\partial W(y)}}{{\partial y}}} \right]}^{ - 1}}{\bm E_{n1}}{\rm{~~~~~~~~~~~~~~~~~~~~~~~~~~~~~~~~~~~~,}}\bm x \to {0^ - }}
\end{array}} \right.
\end{array}\label{eq:6}\end{equation}
where ${\bm E_{n1}} = {\left[ {\begin{array}{*{20}{c}}
1&1& \cdots &1
\end{array}} \right]^T} \in {R^n}$ . It is clear that $\dot {\bm x} $ will be not zero when $\bm x \to \bm 0$ by the theorem of equivalent infinitesimal. Hence, the system (5) may be not stable at 0 as $t \to \infty $.

\textbf{Example 1:} Consider a simple on-order system: $\dot {\bm x}-\dot {\bm x_d}= \bm u$ in which $\bm u=- \frac{1}{{{T_c}}}{\left\| {\frac{{\partial W(\bm x)}}{{\partial \bm  x}}} \right\|^{ - 2}}{\left[ {\frac{{\partial W(\bm x)}}{{\partial \bm x}}} \right]^T}$ and $W(x)=1-exp(\left\|x^m\right\|)$. As shown in Figs. 1 and 2, the simulation is carried out to show the different convergence performance of predefined-time stable system under different parameters $m$. It is clear that the system can converge to the equilibrium point more smoothly in $0<m \le 0.5$ than $0.5<m \le 1$ after $T_c=0.1$. In addition, there is a very interesting performance of system: by consider a traditional idea (m=1 may be an appropriate value in Lemma 1), the chattering degree of controller should decrease as $m$ decreases. However, according to the curves of simulation, the chattering degree of controller in $m=1$ is smaller than that in $0.8 \le m \le 0.9$. Not only that, but this curve is a regular triangular curve, not as irregular as other curves. The mathematical principles of this phenomenon may be explained by \eqref{eq:6}.

\begin{figure}
\begin{center}
\includegraphics[height=5cm]{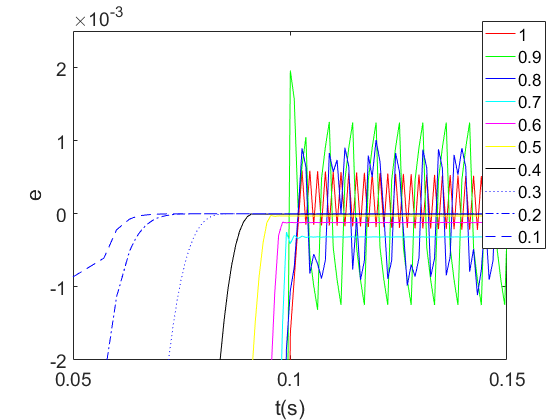}    
\caption{The curve of $(x|\dot x =  - \frac{1}{{{T_c}}}{exp({{{\left| x \right|}^m}}})\frac{x}{{{{\left| x \right|}^m}}})$:}
{$T_c=0.1$, $m=1,0.9,\cdots 0.1$}  
\label{fig1}                                 
\end{center}                                 
\end{figure}

\begin{figure}
\begin{center}
\includegraphics[height=5cm]{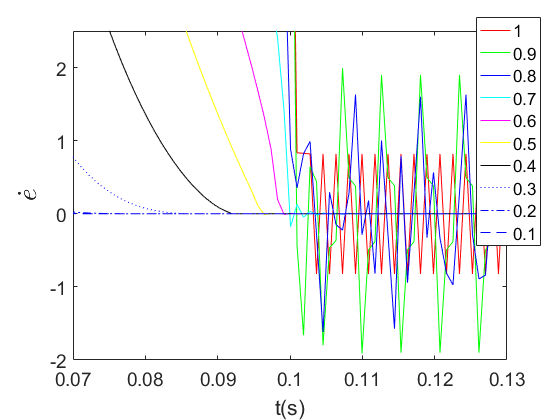}    
\caption{The curve of $(\dot x|\dot x =  - \frac{1}{{{T_c}}}{exp({{{\left| x \right|}^m}}})\frac{x}{{{{\left| x \right|}^m}}})$:}
{$T_c=0.1$, $m=1,0.9,\cdots 0.1$}  
\label{fig2}                                 
\end{center}                                 
\end{figure}

\textbf{Remark 1:} Even by considering the following discontinuous predefined-time stabile system which is given as the case 2 of Definition 3.1 in \cite{sanchez2015predefined}:
\begin{equation}\dot {\bm e} =  - \frac{1}{{{T_c}}}{exp({\left\| {\bm e} \right\|})}{\mathop{\rm sign}\nolimits} ({\bm e}),\label{eq:101}\end{equation}
However, the system is still difficult to converge to the equilibrium point within a predetermined time, because:
\begin{equation}\left\{ {\begin{array}{*{20}{l}}
{\mathop {\lim }\limits_{{\bm e} \to {\bm 0^ + }}  =  - \frac{1}{{{T_c}}}\exp ({{\left\| {\bm e} \right\|}})sign({\bm e}) =  - \frac{1}{{{T_c}}}\exp ({{\left\| {\bm e} \right\|}})}\\
{\mathop {\lim }\limits_{{\bm e} \to {\bm 0^ - }}  =  - \frac{1}{{{T_c}}}\exp ({{\left\| {\bm e} \right\|}})sign({\bm e}) = \frac{1}{{{T_c}}}\exp ({{\left\| {\bm e} \right\|}})}
\end{array}} \right.\label{eq:102}\end{equation}

2)	By considering $0<m<1$, $W({\left\| \bm x \right\|^m})$ still satisfies the condition in Lemma1, and a limit will be exemplified in the following to emphasize $m \in (0,1)$ and $m\ne 1$:
\begin{equation}\begin{array}{l}
\mathop {\lim }\limits_{\bm x \to 0} \dot {\bm x} =  - \frac{1}{{{T_c}}}{\left[ {\frac{{\partial W({{\left\| \bm x \right\|}^m})}}{{\partial {\left\| \bm x \right\|}^m}}\frac{{m{{\left\| \bm x \right\|}^{m - 1}}{x^T}}}{{\left\| \bm x \right\|}}} \right]^{ - 1}}\\
{\rm{  ~~~~~    }}\mathop  = \limits^{y = {{\left\| \bm x \right\|}^m}}  - \frac{1}{{m{T_c}}}{\left[ {\frac{{\partial W(y)}}{{\partial y}}} \right]^{ - 1}}\frac{\bm x}{{{{\left\|\bm  x \right\|}^m}}}\\
{\rm{   ~~~~~~~~    }} = \bm 0
\end{array}.\label{eq:7}\end{equation}

It is clear that $\dot {\bm x}$ will be zero when $\bm x \to \bm 0$ according to the theorem of higher-order infinitesimal. Hence, the system (5) can be stable at 0 as $t \to T_c$.

Consider a candidate Lyapunov function as the following:
\begin{equation}V_1 = \left\| \bm x \right\|.\label{eq:8}\end{equation}

The time derivative of Lyapunov function \eqref{eq:8} can be deduced as:
\begin{equation}\begin{array}{l}
\dot V_1 = \frac{{{\bm x^T}\dot {\bm x}}}{{\left\| \bm x \right\|}}\\
{\rm{ ~~~ }} =  - \frac{1}{{{T_c}}}\frac{{{\bm x^T}}}{{\left\| \bm x \right\|}}{\left[ {\frac{{\partial W({{\left\| \bm x \right\|}^m})}}{{\partial {\left\| \bm x \right\|}^m}}\frac{{m{{\left\| \bm x \right\|}^{m - 1}}{\bm x^T}}}{{\left\| \bm x \right\|}}} \right]^{ - 1}}\\
\mathop  = \limits^{y = {{\left\| \bm x \right\|}^m}}  - \frac{1}{{m{T_c}}}{\left[ {\frac{{\partial W(y)}}{{\partial y}}} \right]^{ - 1}}{\left\| \bm x \right\|^{1 - m}}
\end{array}\label{eq:9}\end{equation}

According to \eqref{eq:3}, the system \eqref{eq:9} may not be stable. The \eqref{eq:3}-(v) may only be a necessary but not sufficient condition, for example:
\begin{equation}\dot {\bm x} =  - \frac{1}{{{mT_c}}}{exp({\left\| {\bm x} \right\|^{m}})}{\left\| {\bm x} \right\|^{1-m}},\label{eq:301}\end{equation}

Although all the conditions is satisfied, but (12) is not stable in $\bm x<0$, hence, \eqref{eq:3}-(v) should be modified as a necessary and sufficient condition: ${\frac{{\partial W({{\left\| \bm x \right\|}^m})}}{{\partial {{\bm x }}}}}>0$, then, (11) can be deduced as:
\begin{equation}\begin{array}{l}
\dot V_1 =  - \frac{1}{{m{T_c}}}{\left[ {\frac{{\partial W({V_1^m})}}{{\partial {V_1^m}}}} \right]^{ - 1}}{\left\|\bm  x \right\|^{1 - m}}\\
~~~ =  - \frac{1}{{m{T_c}}}{\left[ {\frac{{\partial W({V_1^m})}}{{\partial {V_1^m}}}} \right]^{ - 1}}{V_1^{1 - m}}\\
~~~ < 0
\end{array}\label{eq:10}\end{equation}

Then, the settling time can be deduced as:
\begin{equation}\begin{array}{l}
T(\bm x_0) = \int\limits_{V_1({\bm x_0})}^{V_1(\bm 0)} {\frac{1}{{ - \frac{1}{{m{T_c}}}{{\left[ {\frac{{\partial W({V_1^m})}}{{\partial {V_1^m}}}} \right]}^{ - 1}}{V_1^{1 - m}}}}dV_1} \\
 ~~~~~~~~~= {T_c}\int\limits_{V_1{{(\bm 0)}^m}}^{V_1{{({\bm x_0})}^m}} {\frac{{\partial W({V_1^m})}}{{\partial {V_1^m}}}d{V_1^m}} \\
 ~~~~~~~~~= {T_c}\left[ {W(V_1{{({\bm x_0})}^m}) - W(V_1{{(\bm 0)}^m})} \right]
\end{array}.\label{eq:11}\end{equation}
Hence, the max settling time can be direct defined as $T_c$ without any math method of enlarging and reducing. In addition, in my opinion, it is also the essential difference between predefined time stability and fixed time one.

The proof is complete.

\textbf{Definition 2:} Until now, the definition of singularity for terminal sliding mode has been complex and confusing. Hence, we suggest that the singularity for fixed-time terminal sliding mode control can be defined: control input will not be infinite whenever system state is not to infinite. This may be the more proper description of singularity.

\section{Main Results}

\subsection{A Second-order System}

Consider a second-order MIMO system as the following:
\begin{equation}\left\{ {\begin{array}{*{20}{l}}
{{{\dot {\bm x}}_1} = {\bm x_2}}\\
{{{\dot {\bm x}}_2} = f(\bm x,t) + b(\bm x,t)\bm u}\\
{\bm y_o = {\bm x_1}}
\end{array}} \right..\label{eq:12}\end{equation}
where ${\bm x_1} = {\left[ {\begin{array}{*{20}{c}}
{{x_{1,1}}}&{{x_{1,2}}}& \cdots &{{x_{1,n}}}
\end{array}} \right]^T}$, $\bm x = {\left[ {\begin{array}{*{20}{c}}
{{\bm x_1}^T}&{{\bm x_2}^T}
\end{array}} \right]^T}$. $f(\bm x,t):D \to {R^n}$ is a continuous function. $b(\bm x,t):D \to {R^{n \times n}}$ is a continuous function. $\bm u \in {R^n}$ is the control input, $\bm y_o \in {R^n}$ is system output.

\textbf{Assumption 1 (for finite control input):} $\frac{{\partial {{\left[ {\frac{{\partial W(y)}}{{\partial y}}} \right]}^{ - 1}}}}{{\partial y}}>0$ and $\frac{{\partial^2 {{\left[ {\frac{{\partial W(y)}}{{\partial y}}} \right]}^{ - 1}}}}{{\partial y^2}}>0$.

\textbf{Corollary 1:} According to Lemma 2, a predefined-time sliding mode can be designed as:
\begin{equation}\bm s = \dot {\bm e} + \frac{1}{{{T_c}}}{\left\| {\frac{{\partial W({{\left\| \bm e \right\|}^m})}}{{\partial\bm e}}} \right\|^{ - 2}}{\left[ {\frac{{\partial W({{\left\| \bm e \right\|}^m})}}{{\partial \bm e}}} \right]^T}.\label{eq:13}\end{equation}
where $\bm e=\bm y_o-\bm y_d$ in which $\bm y_d$ is the desired system output. Then, a simple sliding mode switching method can be designed to make sliding mode control input nonsingular as:
\begin{equation}m = \left\{ {\begin{array}{*{20}{c}}
{1 + (m - 1)sign(sat(\frac{{\left\|\bm  e \right\|}}{\varepsilon })){\rm{~,~~~ }}\left\| \bm s \right\| \ne 0}\\
{{\rm{       ~~~~~~~~~~~~~~~~~           }}m{\rm{       ~~~~~~~~~~~~~~~~~~~,~~~         }}\left\| \bm s \right\| = 0}
\end{array}} \right..\label{eq:14}\end{equation}
with
\begin{equation}sat(\frac{{\left\| \bm e \right\|}}{\varepsilon }) = \left\{ {\begin{array}{*{20}{c}}
{1{\rm{ ~,~~ }}\left\| \bm e \right\| > \varepsilon }\\
{{\rm{0 ~,~~}}\left\| \bm e \right\| \le \varepsilon }
\end{array}} \right..\label{eq:15}\end{equation}
where $0 < m \le 0.5$ and $\varepsilon $ is a small parameter representing the risk of singularity.
As a direct consequence, the singularity of sliding mode can be overcome by \eqref{eq:14} and \eqref{eq:15} if Assumption 1 is satisfied for $y={\left\| \bm e \right\|}^m$.

Proof:
This subsection mainly discusses the finiteness of control input, hence $b(\bm x,t)$ can set to a n-order unit matrix, $f(\bm x,t)$ can be set to $\bm 0$.

The dynamic of sliding mode surface $\bm s$ can be expressed as:
\begin{equation}\dot {\bm s} = \bm u + \frac{1}{{{T_c}}}\frac{{\partial \left( {{{\left[ {\frac{{\partial W({{\left\| \bm e \right\|}^m})}}{{\partial \bm e}}} \right]}^{ - 1}}} \right)}}{{\partial t}}.\label{eq:16}\end{equation}

1) By considering an arbitrary reaching law $\dot {\bm s} = f(\bm s)$ which satisfies Lemma 2 , then, the control input $\bm u$ can be deduced as:
\begin{equation}\begin{array}{l}
\bm u\mathop  = \limits^{y = {{\left\| \bm e \right\|}^m}} f(\bm s) - \frac{1}{{m{T_c}}}\frac{{\partial \left( {{{\left[ {\frac{{\partial W(y)}}{{\partial y}}} \right]}^{ - 1}}\frac{\bm e}{{{{\left\| \bm e \right\|}^m}}}} \right)}}{{\partial t}}\\
{\rm{  ~~~~~~  }} = f(\bm s)  - \frac{1}{{m{T_c}}}(\frac{{\partial {{\left[ {\frac{{\partial W(y)}}{{\partial y}}} \right]}^{ - 1}}}}{{\partial y}}\frac{{m{{\left\| {\bm e} \right\|}^{m - 1}}{{\bm e}^T}\dot {\bm e}}}{{\left\| {\bm e} \right\|}}\frac{\bm e}{\left\| {\bm e} \right\|^m}\\
~~~~~~~~~~~~~~~~~~~~~~~~~~~~~~~~~~~~~~~~~~~~~~~~ + {\left[ {\frac{{\partial W(y)}}{{\partial y}}} \right]^{ - 1}}\frac{{(1 - m)\dot {\bm e}}}{{{{\left\| {\bm e} \right\|}^m}}})\\
{\rm{  ~~~~~~  }} = f(\bm s) - \frac{1}{{{T_c}}}\frac{{\partial {{\left[ {\frac{{\partial W(y)}}{{\partial y}}} \right]}^{ - 1}}}}{{\partial y}}\dot {\bm e} - \frac{{(1 - m)}}{{m{T_c}}}{\left[ {\frac{{\partial W(y)}}{{\partial y}}} \right]^{ - 1}}\frac{{\dot {\bm e}}}{{{{\left\| \bm e \right\|}^m}}}
\end{array}\label{eq:17}\end{equation}

Then, the singularity will come (generally, $\bm {\dot e} \to \infty$ will not be considered):
\begin{equation}\begin{array}{l}
\mathop {\lim }\limits_{\bm e \to 0,\dot {\bm e} \ne 0} \bm u\mathop  = \limits^{y = {{\left\| \bm e \right\|}^m}} f(\bm s) -\frac{1}{{{T_c}}}\frac{{\partial {{\left[ {\frac{{\partial W(y)}}{{\partial y}}} \right]}^{ - 1}}}}{{\partial y}}\dot {\bm e} - \\
~~~~~~~~~~~~~~~~~~~~~~~~~~~~~~~~~~~~~~~~~~~~~~~\frac{{(1 - m)}}{{m{T_c}}}{\left[ {\frac{{\partial W(y)}}{{\partial y}}} \right]^{ - 1}}\frac{{\dot {\bm e}}}{{{{\left\| \bm 0 \right\|}^m}}}\\
~~~~~~~~~~~~~~~~~ = f(\bm s) - \frac{1}{{{T_c}}}\frac{{\partial {{\left[ {\frac{{\partial W(y)}}{{\partial y}}} \right]}^{ - 1}}}}{{\partial y}}\dot {\bm e} - \infty {\mathop{\rm sign}\nolimits} (\dot {\bm e})\\
~~~~~~~~~~~~~~~~~ = f(\bm s) - \infty {\mathop{\rm sign}\nolimits} (\dot {\bm e})
\end{array}\label{eq:18}\end{equation}

Let me break off the proof for a moment, an interesting phenomenon can be discussed to express my alternative understanding of singularity.

a). By considering $f(\bm s)$ is finite, then,
\begin{equation}\mathop {\lim }\limits_{\bm e \to 0,\dot {\bm e} \ne 0}\bm u =  - \infty {\mathop{\rm sign}\nolimits} (\dot {\bm e}).\label{eq:19}\end{equation}

b). By considering $\left\|f(\bm s)\right\|$ is the inequivalent infinity of $\left\| \frac{1}{{{T_c}}}\frac{{\partial {{\left[ {\frac{{\partial W(y)}}{{\partial y}}} \right]}^{ - 1}}}}{{\partial y}}\dot {\bm e} + \infty {\mathop{\rm sign}\nolimits} (\dot {\bm e})\right\|$ in \eqref{eq:18} and ${\Psi _i} = {\mathop{\rm sign}\nolimits} (f({s_i})) - {\mathop{\rm sign}\nolimits} ({\dot e_i})$, then,
\begin{equation}
\mathop {\lim }\limits_{\bm e \to 0,\dot {\bm e} \ne 0} \bm u = \pm \infty
\label{eq:302}\end{equation}

c). By considering $\left\| f(\bm s)\right\|$ is the equivalent infinity of $ \left\| \frac{1}{{{T_c}}}\frac{{\partial {{\left[ {\frac{{\partial W(y)}}{{\partial y}}} \right]}^{ - 1}}}}{{\partial y}}\dot {\bm e} + \infty {\mathop{\rm sign}\nolimits} (\dot {\bm e})\right\|$ in \eqref{eq:18} and ${\Psi _i} = {\mathop{\rm sign}\nolimits} (f({s_i})) - {\mathop{\rm sign}\nolimits} ({\dot e_i})$, then,
\begin{equation}\begin{array}{l}
\mathop {\lim }\limits_{\bm e \to 0,\dot {\bm e} \ne 0} \bm u = \infty \left[ {{\mathop{\rm sign}\nolimits} (f(\bm s)) - {\mathop{\rm sign}\nolimits} (\dot {\bm e})} \right]\\
~~~~~~~~~~~~~~ = {\left[ {\begin{array}{*{20}{c}}
{\infty {\Psi _1}}&{\infty {\Psi _2}}& \cdots &{\infty {\Psi _n}}
\end{array}} \right]^T}\\
\end{array}\label{eq:20}\end{equation}
In case of (c), it is clear that the control input $\bm u$ may be zero if each $f({s_i}){\dot e_i} > 0$.

\textbf{Remark 2:} In the perspective of control engineering, the infinite initial condition may not be controlled, however, according to \eqref{eq:18} and \eqref{eq:20}, the control input may be finite in a very strict condition as the equivalent infinity in (c). In addition, according to $\dot {\bm s} = f(\bm s)$ in \eqref{eq:17}, if each ${\dot s_i}{\dot x_i} > 0$ and there is no unknown disturbance, the strict condition may be relaxed as: $\mathop {\lim }\limits_{{\bm e} \to \frac{1}{\bm c},\dot {\bm e}{\rm{ = }}\bm c} \left\| \bm u \right\|\mathop {\mathop  = \limits_{f({s_i}){{\dot e}_i} > 0} }\limits^{y = {{\left\| {\bm e} \right\|}^m}} \left\| {f(\bm s)} \right\| - (\left\| {\frac{1}{{{T_c}}}\frac{{\partial {{\left[ {\frac{{\partial W(y)}}{{\partial y}}} \right]}^{ - 1}}}}{{\partial y}}\dot {\bm e}} \right\| + \left\| {\frac{{(1 - m)}}{{m{T_c}}}{{\left[ {\frac{{\partial W(y)}}{{\partial y}}} \right]}^{ - 1}}\frac{{\dot {\bm e}}}{{{{\left\| 0 \right\|}^m}}}} \right\|) \le {\bm{~a ~finite~ constant~}}$ in which $\left\| \bm c \right\|$ is a big number. i.e. $\left\| \bm c \right\| = 3.0 \times {10^8}$. Hence, a special class of infinite initial condition may be controlled to be stable by a global finite control input if controller can make ${\dot s_i}{\dot e_i} > 0$ all the time. The Remark 2 is only a conjecture, but it also reveals some properties of singularity.

Go on the proof. Then, according to \eqref{eq:18}, \eqref{eq:19}, \eqref{eq:302} and \eqref{eq:20}, we know that the possible reason of infinite control inputs. By considering actual engineer, $f(\bm s)$ and $\dot {\bm e}$ is finite, hence, we set $m=1$ when $\bm e \to 0$ and $\bm s \ne 0$:
\begin{equation}\begin{array}{l}\mathop {\lim }\limits_{{\bm e} \to 0,\dot {\bm e} \ne 0} \left\| {\bm u} \right\|\mathop {\mathop  = \limits_{m = 1} }\limits^{y = {{\left\| {\bm e} \right\|}^m}} f(\bm s) - \frac{1}{{{T_c}}}\frac{{\partial {{\left[ {\frac{{\partial W(y)}}{{\partial y}}} \right]}^{ - 1}}}}{{\partial y}}\dot {\bm e}\\
~~~~~~~~~~~~~~~~~~~~~ = {\bm{~a ~finite ~constant}}
\end{array}\label{eq:22}\end{equation}

Then, by considering a predefined-time reaching law $f(\bm s) =  - \frac{1}{{{T_c}}}{\left\| {\frac{{\partial W({{\left\|\bm  s \right\|}^m})}}{{\partial {\left\|\bm  s \right\|}^m}}} \right\|^{ - 2}}{\left[ {\frac{{\partial W({{\left\| \bm s \right\|}^m})}}{{\partial \bm s}}} \right]^T}$, the candidate Lyapunov function is $V_s=\left\|\bm s\right\|$ and the sliding mode existence condition can be expressed:
\begin{equation}\begin{array}{l}
{{\dot V}_s} = \frac{{{{\bm s}^T} \bm {\dot s}}}{{\left\| {\bm s} \right\|}}\\
{\rm{ ~~~ }} =  - \frac{1}{{{T_c}}}\frac{{{{\bm s}^T}}}{{\left\| {\bm s} \right\|}}{\left[ {\frac{{\partial W({{\left\| {\bm s} \right\|}^m})}}{{\partial {{\left\|\bm  s \right\|}^m}}}\frac{{m{{\left\| {\bm s} \right\|}^{m - 1}}{\bm s^T}}}{{\left\| {\bm s} \right\|}}} \right]^{ - 1}}\\
 ~~~=  - \frac{1}{{m{T_c}}}{\left[ {\frac{{\partial W({V_s}^m)}}{{\partial {V_s}^m}}} \right]^{ - 1}}{V_s}^{1 - m}\\
 ~~~< 0
\end{array}\label{eq:23}\end{equation}
The reaching time can be deduced by:
\begin{equation}\begin{array}{l}
T = \int\limits_{{V_s}({\bm e_0})}^{{V_s}({\bm e_{s = 0}})} {\frac{1}{{ - \frac{1}{{m{T_c}}}{{\left[ {\frac{{\partial W({V_s}^m)}}{{\partial {V_s}^m}}} \right]}^{ - 1}}{V_s}^{1 - m}}}d{V_s}} \\
 ~~= {T_c}\int\limits_{{V_s}{{({\bm e_{s = 0}})}^m}}^{{V_s}{{({\bm e_0})}^m}} {\frac{{\partial W({V^m})}}{{\partial {V^m}}}d{V_s}^m} \\
 ~~= {T_c}\left[ {W({V_s}{{({{\bm e}_0})}^m}) - W({V_s}{{({{\bm e}_{s = 0}})}^m})} \right]
\end{array}.\label{eq:24}\end{equation}
Hence, the minimum upper bound of reaching time is $T_c$.

\textbf{Result 1:} the parameter $m$ satisfied the condition of \eqref{eq:14} will not affect the minimum upper bound of reaching time.

\begin{figure}
\begin{center}
\includegraphics[height=5cm]{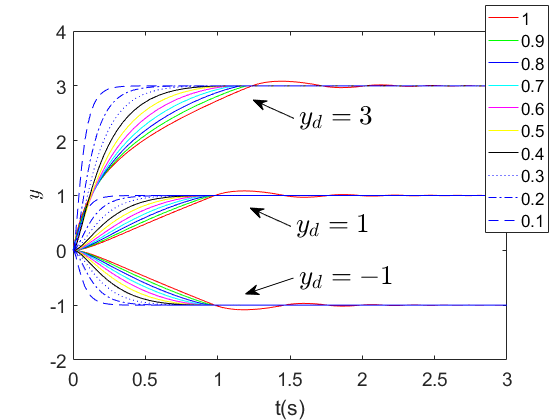}    
\caption{The curve of $\bm y$ in which $W(\left\|\bm x\right\|^m)=1-exp(\left\|\bm x\right\|^m)$:}
{$T_c=1$, $m=1,0.9,\cdots 0.1$}  
\label{fig1}                                 
\end{center}                                 
\end{figure}

2)	When sliding mode surface reaches 0, then, \eqref{eq:13} can be written as:
\begin{equation}{\bm s}=\bm 0=\dot {\bm e} + \frac{1}{{m{T_c}}}{\left\| {\frac{{\partial W({{\left\| {\bm e} \right\|}^m})}}{{\partial {\bm e}}}} \right\|^{ - 2}}{\left[ {\frac{{\partial W({{\left\| {\bm e} \right\|}^m})}}{{\partial {\bm e}}}} \right]^T}.\label{eq:25}\end{equation}

Then control input could be expressed by \eqref{eq:17} and \eqref{eq:25} as the following:
\begin{equation}\begin{array}{l}
\bm u\mathop  = \limits^{y = {{\left\| {\bm e} \right\|}^m}} f(\bm s) - \frac{1}{{m{T_c}}}(\frac{{\partial {{\left[ {\frac{{\partial W(y)}}{{\partial y}}} \right]}^{ - 1}}}}{{\partial y}}m\dot {\bm e} \\
~~~~~~~~~+ {\left[ {\frac{{\partial W(y)}}{{\partial y}}} \right]^{ - 1}}\frac{{(m - 1)\frac{1}{{m{T_c}}}{{\left\| {\frac{{\partial W({{\left\| {\bm e} \right\|}^m})}}{{\partial {{\left\| {\bm e} \right\|}^m}}}} \right\|}^{ - 2}}{{\left[ {\frac{{\partial W({{\left\| {\bm e} \right\|}^m})}}{{\partial {\bm e}}}} \right]}^T}}}{{{{\left\| {\bm e} \right\|}^m}}})\\
~~~~~~ = f(\bm s) - \frac{1}{{{T_c}}}\frac{{\partial {{\left[ {\frac{{\partial W(y)}}{{\partial y}}} \right]}^{ - 1}}}}{{\partial y}}\dot {\bm e} - \frac{{(m - 1)}}{{{m^2}{T^2}_c}}{\left[ {\frac{{\partial W(y)}}{{\partial y}}} \right]^{ - 2}}\frac{{\bm e}}{{{{\left\| {\bm e} \right\|}^{2m}}}}
\end{array}\label{eq:26}\end{equation}

According to higher-order infinitesimal theorem, it is clear that control input $\bm u$ in \eqref{eq:26} will not be infinite in any actual conditions when $0 < m \le \frac{1}{2}$.

Then a candidate Lyapunov function can be shown as:
\begin{equation}V_2 = \left\| {\bm e} \right\|.\label{eq:27}\end{equation}

The time derivative of Lyapunov function \eqref{eq:27} can be deduced as:
\begin{equation}\begin{array}{l}
\dot V_2 = \frac{{{{\bm e}^T}\dot {\bm e}}}{{\left\| {\bm e} \right\|}}\\
{\rm{ ~~~ }} =  - \frac{1}{{{T_c}}}\frac{{{{\bm e}^T}}}{{\left\| {\bm e} \right\|}}{\left[ {\frac{{\partial W({{\left\| {\bm e} \right\|}^m})}}{{\partial {{\left\| {\bm e} \right\|}^m}}}\frac{{m{{\left\| {\bm e} \right\|}^{m - 1}}{{\bm e}^T}}}{{\left\| {\bm e} \right\|}}} \right]^{ - 1}}\\
{\rm{ }}\mathop  = \limits^{y = {{\left\| {\bm e} \right\|}^m}}  - \frac{1}{{m{T_c}}}{\left[ {\frac{{\partial W(y)}}{{\partial y}}} \right]^{ - 1}}{\left\| {\bm e} \right\|^{1 - m}}
\end{array}\label{eq:28}\end{equation}

When $m=1$:
\begin{equation}\begin{array}{l}
\dot V_2 =  - \frac{1}{{{T_c}}}{\left[ {\frac{{\partial W(V_2)}}{{\partial V_2}}} \right]^{ - 1}}\\
~~~ < 0
\end{array}.\label{eq:29}\end{equation}

The settling time function of \eqref{eq:29} can be deduced:
\begin{equation}\begin{array}{l}
T = \int\limits_{V_2({{\bm e}_0})}^{V_2(\bm 0)} {\frac{1}{{ - \frac{1}{{{T_c}}}{{\left[ {\frac{{\partial W(V_2)}}{{\partial V_2}}} \right]}^{ - 1}}}}dV_2} \\
 ~~= {T_c}\int\limits_{V_2(\bm 0)}^{V_2({{\bm e}_0})} {\frac{{\partial W(V_2)}}{{\partial V_2}}dV_2} \\
 ~~= {T_c}\left[ {W(V_2({{\bm e}_0})) - W(V_2(\bm 0))} \right]
\end{array}.\label{eq:30}\end{equation}

Note that $m \ne 1$ when $s=0$, because $\dot V \ne 0$ if $V=0$. In addition, by considering $f(\bm x,t)=-\bm x_2$ and $b(\bm x,t)=\bm E_{3\times 3}$. The output of system with the controller designed by \eqref{eq:17} can be showed in Fig. 3. According to Fig. 3, it is clear that the system errors will not converge to zero within $2T_c$ only if $m=1$, and it converge wavily to zero at about 2.5s.

When $0<m \le \frac{1}{2}$.
\begin{equation}\begin{array}{l}
\dot V_2 =  - \frac{1}{{m{T_c}}}{\left[ {\frac{{\partial W({V_2^m})}}{{\partial {V_2^m}}}} \right]^{ - 1}}{V_2^{1 - m}}\\
 ~~~< 0
\end{array}.\label{eq:31}\end{equation}

The settling time function of \eqref{eq:31} can be deduced:
\begin{equation}\begin{array}{l}
T = \int\limits_{V_2({{\bm e}_0})}^{V_2(\bm 0)} {\frac{1}{{ - \frac{1}{{m{T_c}}}{{\left[ {\frac{{\partial W({V_2^m})}}{{\partial {V_2^m}}}} \right]}^{ - 1}}{V_2^{1 - m}}}}dV_2} \\
~~ = {T_c}\int\limits_{V_2{{(\bm 0)}^m}}^{V_2{{({{\bm e}_0})}^m}} {\frac{{\partial W({V_2^m})}}{{\partial {V_2^m}}}d{V_2^m}} \\
~~ = {T_c}\left[ {W(V_2{{({{\bm e}_0})}^m}) - W(V_2{{(\bm 0)}^m})} \right]
\end{array}.\label{eq:32}\end{equation}

The upper bounds of reaching time in the case of \eqref{eq:29} and \eqref{eq:31} are same. However, these equilibrium states may be different: \eqref{eq:29} is  ${(V_2,\dot V_2)_{V_2 = 0}} = {(V_2, - \frac{1}{{{T_c}}}{\left[ {\frac{{\partial W(V_2)}}{{\partial V_2}}} \right]^{ - 1}})_{V_2 = 0}} = (0, - \frac{1}{{{T_c}}}{\left[ {\frac{{\partial W(V_2 \to {0^{\rm{ + }}})}}{{\partial (V_2 \to {0^{\rm{ + }}})}}} \right]^{ - 1}})$, \eqref{eq:31} is ${(V_2,\dot V_2)_{V_2 = 0}} ~~~= {(V_2, - \frac{1}{{{T_c}}}{\left[ {\frac{{\partial W(V_2)}}{{\partial V_2}}} \right]^{ - 1}}{V_2^{1 - m}})_{V_2 = 0}} = (0,0)$. In addition, none of the $y$ derivative of the known predefined-time stability functions is 0 if $V_2 \to {0}$. Hence, the parameter should be $0<m \le \frac{1}{2}$ when $\left\| \bm s \right\| = 0$.

In order to make control system nonsingular and predefined-time stable, the parameter $m$ should be changed as \eqref{eq:14}.

The proof is complete.

\textbf{Remark 3:} By considering a design method of semi-global predefined-time system in \cite{jimenez2017semi,jimenez2018semi,sanchez2018class}. A class of candidate Lyapunov function can be set as:
\begin{equation}V =  {\left\| {\bm e} \right\|^m}.\label{eq:33}\end{equation}
Then, the predefined-time system could be designed as:
\begin{equation}\dot {\bm e} =  - \frac{1}{T_c}{exp({\left\| {\bm e} \right\|^{mp}})}\frac{\bm e}{\left\| {\bm e} \right\|^{mp}}.\label{eq:34}\end{equation}
In this case, the parameters $m$ and $p$ will be considered combined into a single parameter $m=mp$.

\textbf{Remark 4:} In actual control engineering using the SMC control, $\bm s$ will not be zero in the presence of perturbations. Hence, the sliding mode switching strategy in \eqref{eq:14} will work frequently, and then, an unavoidable chattering of control input will occur when sliding mode switches. We suggest that sliding mode switching strategy in \eqref{eq:14} will work only if control input $\bm u$ is real large:
\begin{equation}\varepsilon  = \left\{ {\begin{array}{*{20}{c}}
{\varepsilon ~,~~~\left\| \bm u \right\| > {u_{\max }}}\\
{0~,~~~~others~~~~~~~}
\end{array}} \right..\label{eq:35}\end{equation}
where $u_{max}$ is the maximum allowable value of control input.

\subsection{A N-order System}

Consider a n-order MIMO system as the following:
\begin{equation}\left\{ {\begin{array}{*{20}{l}}
{{{\dot {\bm x}}_i} = {\bm x_{i+1}}}~;~~~i=1,2 \cdots,n-1 \\
{{{\dot {\bm x}}_n} = f(\bm x,t) + b(\bm x,t)(\bm u+\bm u_r) + \bm d(t)}\\
{\bm y_o = {\bm x_1}}
\end{array}} \right.\label{eq:36}\end{equation}
where ${\bm x} = {\left[ {\begin{array}{*{20}{c}}
{{\bm x_{1}^T}}&{{\bm x_{2}^T}}& \cdots &{{\bm x_{n}^T}}
\end{array}} \right]^T}$, $\bm u_r \in R^n$ is a robust controller, and an accurate estimation ${\hat {\bm d}}(t)$ of $\bm d(t)$ will be considered to design $\bm u_r \in R^n=b^{-1}(\bm x,t)\left\|{\hat {\bm d}}(t)\right\|sign(\bm s_{n-1})$.

Consider a tracking problem $\bm y \to {\bm y_d}$, and the error vector can be defined as:
\begin{equation}\bm \Im = \bm x - {\bm x_d} = {\left[ {\begin{array}{*{20}{c}}
{{\bm e^T}}&{{{\dot {\bm e}}^T}}& \cdots &{{\bm e^{(n - 1)}}^T}
\end{array}} \right]^T}.\label{eq:37}\end{equation}

\textbf{Corollary 2:} The singularity of n-order system can be described as: control input $\bm u$ will be infinite if $\bm e \to 0$ and $\left\| \bm \Im \right\| \ne 0$. A nonsingular recursive predefined-time sliding mode can be designed by Lemma 2 as:
\begin{equation}\begin{array}{l}
~~~~{{\bm s}_1} = {{\dot {\bm s}}_0} + \frac{1}{{{T_c}}}{\left\| {\frac{{\partial W({{\left\| {{{\bm s}_0}} \right\|}^{m_0}})}}{{\partial {{\bm s}_0}}}} \right\|^{ - 2}}{\left[ {\frac{{\partial W({{\left\| {{{\bm s}_0}} \right\|}^{m_0}})}}{{\partial {{\bm s}_0}}}} \right]^T}\\
~~~~{{\bm s}_2} = {{\dot {\bm s}}_1} + \frac{1}{{{T_c}}}{\left\| {\frac{{\partial W({{\left\| {{{\bm s}_1}} \right\|}^{m_1}})}}{{\partial {{\bm s}_1}}}} \right\|^{ - 2}}{\left[ {\frac{{\partial W({{\left\| {{{\bm s}_1}} \right\|}^{m_1}})}}{{\partial {{\bm s}_1}}}} \right]^T}\\
{\rm{            ~~~~~~~~~~~~~~~~~~~~~~~~~~~~~~~~~~~~              }} \vdots \\
{{\bm s}_{n - 1}} = {{\dot {\bm s}}_{n - 2}} + \frac{1}{{{T_c}}}{\left\| {\frac{{\partial W({{\left\| {{{\bm s}_{n - 2}}} \right\|}^{m_{n-2}}})}}{{\partial {{\bm s}_{n - 2}}}}} \right\|^{ - 2}}{\left[ {\frac{{\partial W({{\left\| {{{\bm s}_{n - 2}}} \right\|}^{m_{n-2}}})}}{{\partial {{\bm s}_{n - 2}}}}} \right]^T}\\
~~~~~\bm 0 = {{\dot {\bm s}}_{n - 1}} + \frac{1}{{{T_c}}}{\left\| {\frac{{\partial W({{\left\| {{{\bm s}_{n - 1}}} \right\|}^{m_{n-1}}})}}{{\partial {{\bm s}_{n - 1}}}}} \right\|^{ - 2}}{\left[ {\frac{{\partial W({{\left\| {{{\bm s}_{n - 1}}} \right\|}^{m_{n-1}}})}}{{\partial {{\bm s}_{n - 1}}}}} \right]^T},
\end{array}\label{eq:38}\end{equation}
in which $\bm s_0=e$ and the parameters satisfy the following condition:
\begin{equation}m_k = \left\{ {\begin{array}{*{20}{c}}
{1 + (m_k - 1)sign(sat(\frac{{\left\|\bm  s_k \right\|}}{\varepsilon })){\rm{~,~~~ }}\left\| \bm s_{k+1} \right\| \ne 0}\\
{{\rm{       ~~~~~~~~~~~~~~~~~~           }}m_k{\rm{       ~~~~~~~~~~~~~~~~~~~~,~~~         }}\left\| \bm s_{k+1} \right\| = 0}
\end{array}} \right..\label{eq:112}\end{equation}
with
\begin{equation}sat(\frac{{\left\| \bm s_k \right\|}}{\varepsilon }) = \left\{ {\begin{array}{*{20}{c}}
{1{\rm{ ~,~~ }}\left\| \bm s_k \right\| > \varepsilon }\\
{{\rm{0 ~,~~}}\left\| \bm s_k \right\| \le \varepsilon }
\end{array}} \right.\label{eq:113}\end{equation}
where $k=0,\cdots, n-2$. $0 < m_{k} \le \frac{1}{n-k}$ and $0<m_{n-1}<1$. $\varepsilon $ is a small parameter representing the risk of singularity which satisfy the Remark 4.

\textbf{Assumption 2:} $\frac{{{\partial ^{(n)}}{{\left[ {\frac{{\partial W({{\left\| {{\bm s_{n - 1}}} \right\|}^m})}}{{\partial {{\left\| {{\bm s_{n - 1}}} \right\|}^m}}}} \right]}^{ - 1}}}}{{\partial {{({{\left\| {{\bm s_{n - 1}}} \right\|}^m})}^{(n)}}}}$ should be bounded in finite area for $0<m\le 1$, $n \in R$ and $\bm s_{n-1} \in R^n\backslash ({\infty })$. i.e. $W(\left\| x \right\|^m)=1-e^{-{\left\| x \right\|^m}}$.

\textbf{Corollary 2 is valid only if Assumption 2 is satisfied.}

Then, according to Corollary 2, a sliding mode reaching law can be chose as:
\begin{equation}{\dot {\bm s}_{n - 1}} =  - \frac{1}{{{T_c}}}{\left\| {\frac{{\partial W({{\left\| {{{\bm s}_{n - 1}}} \right\|}^{m_{n-1}}})}}{{\partial {{\bm s}_{n - 1}}}}} \right\|^{ - 2}}{\left[ {\frac{{\partial W({{\left\| {{{\bm s}_{n - 1}}} \right\|}^{m_{n-1}}})}}{{\partial {{\bm s}_{n - 1}}}}} \right]^T}.\label{eq:39}\end{equation}

The nonsingular predefined-time sliding mode control input $\bm u$ for n-order MIMO system \eqref{eq:36} with global finite control input can be designed by \eqref{eq:26}, \eqref{eq:38} and \eqref{eq:39} as:
\begin{equation}\begin{array}{l}
\bm u =  - {b^{ - 1}}(\bm x,t)\left[ {f(\bm x,t)+y_d^{(n)}} \right.\\
 ~~~~~~+ \sum\limits_{k = 0}^{n - 2} {(\frac{1}{{{T_c}}}{{\left\| {\frac{{\partial W({{\left\| {{\bm s_k}} \right\|}^{m_{k}}})}}{{\partial {\bm s_k}}}} \right\|}^{ - 2}}{{\left[ {\frac{{\partial W({{\left\| {{\bm s_k}} \right\|}^{m_{k}}})}}{{\partial {\bm s_k}}}} \right]}^T}} {)^{(n - k - 1)}}\\
~~~~~~\left. { + \frac{1}{{{T_c}}}{{\left\| {\frac{{\partial W({{\left\| {{{\bm s}_{n - 1}}} \right\|}^{m_{n-1}}})}}{{\partial {{\bm s}_{n - 1}}}}} \right\|}^{ - 2}}{{\left[ {\frac{{\partial W({{\left\| {{{\bm s}_{n - 1}}} \right\|}^{m_{n-1}}})}}{{\partial {{\bm s}_{n - 1}}}}} \right]}^T}} \right]\\
 ~~~=  - {b^{ - 1}}(\bm x,t)\left[ {f(\bm x,t)} \right.\\
 ~~~~~~+ \sum\limits_{k = 0}^{n - 2} {\frac{1}{{{m_{k}}{T_c}}}({{\left[ {\frac{{\partial W({{\left\| {{{\bm s}_k}} \right\|}^{m_{k}}})}}{{\partial {{\left\| {{{\bm s}_k}} \right\|}^{m_{k}}}}}} \right]}^{ - 1}}} \frac{{{{\bm s}_k}}}{{{{\left\| {{{\bm s}_k}} \right\|}^{m_{k}}}}}{)^{(n - k - 1)}}\\
~~~~~~\left. { + \frac{1}{{{m_{n-1}}{T_c}}}{{\left[ {\frac{{\partial W({{\left\| {{{\bm s}_{n - 1}}} \right\|}^{m_{n-1}}})}}{{\partial {{\left\| {{{\bm s}_{n - 1}}} \right\|}^{m_{n-1}}}}}} \right]}^{ - 1}}\frac{{{{\bm s}_{n - 1}}}}{{{{\left\| {{{\bm s}_{n - 1}}} \right\|}^{m_{n-1}}}}}} \right]
\end{array}\label{eq:40}\end{equation}
where $(\cdot )^{(n-k-1)}=\frac{\partial^{(n-k-1)}(\cdot)}{\partial t^{(n-k-1)}}$.

Proof of predefined-time stable control system:

Consider a candidate Lyapunov function as:
\begin{equation}{V_1} = \left\| {{\bm s_{n - 1}}} \right\|\label{eq:103}\end{equation}

The time derivative of Lyapunov function \eqref{eq:103} can be deduced as:
\begin{equation}\begin{array}{l}
{{\dot V}_1} = \frac{{{{\bm s}_{n - 1}}^T{{\dot {\bm s}}_{n - 1}}}}{{\left\| {{{\bm s}_{n - 1}}} \right\|}}\\
~~~ =\frac{{{\bm s_{n - 1}}^T}}{{\left\| {{\bm s_{n - 1}}} \right\|}}\left\{ {f(\bm x,t) + b(\bm x,t)(\bm u + {\bm u_r}) + \bm d(t)} \right.\\
~~~~~~\left. { + \sum\limits_{k = 0}^{n - 2} {\frac{1}{{{m_{k}}{T_c}}}{{({{\left[ {\frac{{\partial W({{\left\| {{\bm s_k}} \right\|}^{m_{k}}})}}{{\partial {{\left\| {{\bm s_k}} \right\|}^{m_{k}}}}}} \right]}^{ - 1}}\frac{{{\bm s_k}}}{{{{\left\| {{\bm s_k}} \right\|}^{m_{k}}}}})}^{(n - k - 1)}}} } \right\}\\
~~~ = \frac{{{{\bm s}_{n - 1}}^T}}{{\left\| {{{\bm s}_{n - 1}}} \right\|}}\left\{ {f(\bm x,t) - b(\bm x,t){b^{ - 1}}(\bm x,t)\left[ {f(\bm x,t) } \right.} \right.\\
~~~~~~ + \sum\limits_{k = 0}^{n - 2} {\frac{1}{{{m_{k}}{T_c}}}{{({{\left[ {\frac{{\partial W({{\left\| {{{\bm s}_k}} \right\|}^{m_{k}}})}}{{\partial {{\left\| {{{\bm s}_k}} \right\|}^{m_{k}}}}}} \right]}^{ - 1}}\frac{{{{\bm s}_k}}}{{{{\left\| {{{\bm s}_k}} \right\|}^{m_{k}}}}})}^{(n - k - 1)}}} \\
~~~~~~ + \frac{1}{{{m_{n-1}}{T_c}}}{{\left[ {\frac{{\partial W({{\left\| {{{\bm s}_{n - 1}}} \right\|}^{m_{n-1}}})}}{{\partial {{\left\| {{{\bm s}_{n - 1}}} \right\|}^{m_{n-1}}}}}} \right]}^{ - 1}}\frac{{{{\bm s}_{n - 1}}}}{{{{\left\| {{{\bm s}_{n - 1}}} \right\|}^{m_{n-1}}}}} \\
~~~~~~ \left. {+ \left\|{\hat  {\bm d}(t)}\right\|{\mathop{\rm sign}\nolimits} ({{\bm s}_{n - 1}})} \right] + \bm d(t)\\
~~~~~~ \left. { + \sum\limits_{k = 0}^{n - 2} {\frac{1}{{{m_{k}}{T_c}}}{{({{\left[ {\frac{{\partial W({{\left\| {{{\bm s}_k}} \right\|}^{m_{k}}})}}{{\partial {{\left\| {{{\bm s}_k}} \right\|}^{m_{k}}}}}} \right]}^{ - 1}}\frac{{{{\bm s}_k}}}{{{{\left\| {{{\bm s}_k}} \right\|}^{m_{k}}}}})}^{(n - k - 1)}}} } \right\}\\
~~~{\rm{  }} =  - \frac{1}{{{m_{n-1}}{T_c}}}\frac{{{{\bm s}_{n - 1}}^T}}{{\left\| {{{\bm s}_{n - 1}}} \right\|}}{\left[ {\frac{{\partial W({{\left\| {{{\bm s}_{n - 1}}} \right\|}^{m_{n-1}}})}}{{\partial {{\left\| {{{\bm s}_{n - 1}}} \right\|}^{m_{n-1}}}}}} \right]^{ - 1}}\frac{{{{\bm s}_{n - 1}}}}{{{{\left\| {{{\bm s}_{n - 1}}} \right\|}^{m_{n-1}}}}}\\
~~~ =  - \frac{1}{{{m_{n-1}}{T_c}}}{\left[ {\frac{{\partial W({{\left\| {\bm s_{n-1}} \right\|}^{m_{n-1}}})}}{{\partial {{\left\| {\bm s_{n-1}} \right\|}^{m_{n-1}}}}}} \right]^{ - 1}}{\left\| {\bm s_{n-1}} \right\|^{1 - {m_{n-1}}}}\\
~~~ =  - \frac{1}{{{m_{n-1}}{T_c}}}{\left[ {\frac{{\partial W({V_1}^{m_{n-1}})}}{{\partial {V_1}^{m_{n-1}}}}} \right]^{ - 1}}{V_1}^{1 - {m_{n-1}}}
\end{array}\label{eq:104}\end{equation}

Then, the reaching time of $\bm s_{n-1}$ can be deduced as:
\begin{equation}\begin{array}{l}
{T_{n - 1}}({\bm s_{n-1}(\bm 0)})\\
 ~~~= \int\limits_{{V_1}({\bm s_{n-1}(0)})}^{{V_1}({{\bm 0}})} {\frac{1}{{ - \frac{1}{{{m_{n-1}}{T_c}}}{{\left[ {\frac{{\partial W({V_1}^{m_{n-1}})}}{{\partial {V_1}^{m_{n-1}}}}} \right]}^{ - 1}}{V_1}^{1 - {m_{n-1}}}}}d} {V_1}\\
~~~= {T_c}\int\limits_{{V_1}({{{} \bm  0}})^{m_{n-1}}}^{{V_1}({\bm s_{n-1}(0)})^{m_{n-1}}} {\frac{{\partial W({V_1}^{m_{n-1}})}}{{\partial {V_1}^{m_{n-1}}}}d} {V_1}^{m_{n-1}}\\
~~~ = {T_c}(W({V_1}{({{\bm s_{n-1}}(\bm 0)})^{m_{n-1}}}) - W({V_1}{({{\bm  0}})^{m_{n-1}}}))
\end{array}\label{eq:105}\end{equation}

It is clear that ${\bm s_{n - 1}} = 0$ after $\sum\limits_{k = 1}^{n-1} {{T_{n - k}}(\bm s_{n-k}(\bm 0))} $. In a similar way of \eqref{eq:104} and \eqref{eq:105}, each reaching time of ${s_{n - k}}$ can be deduced as ${T_c}(W({V_1}{({\bm x_{{s_{n - k}} = 0}})^{m_{n-1}}}) - W({V_1}{({\bm x_{{s_{n - k - 1}} = 0}})^{m_{n-1}}}))$, and $k = 1,2 \cdots n - 1$.
Then, the minimum upper bound of reaching time is $(n - 1){T_c}$.

After $(n - 1){T_c}$, the sliding mode surface ${\bm s_1} = 0$ has been reached, and \eqref{eq:38} can be expressed as:
\begin{equation}\dot {\bm e} =  - \frac{1}{{m_0{T_c}}}{\left\| {\frac{{\partial W({{\left\| {\bm e} \right\|}^{m_0}})}}{{\partial {\bm e}}}} \right\|^{ - 2}}{\left[ {\frac{{\partial W({{\left\| {\bm e}\right\|}^{m_0}})}}{{\partial {\bm e}}}} \right]^T}\label{eq:106}\end{equation}

The time derivative of a candidate Lyapunov function: $V_2={\left\| {\bm e} \right\|}$ can be deduced as:
\begin{equation}\begin{array}{l}
{{\dot V}_2} = \frac{{{{\bm e}^T}\dot {\bm e}}}{{\left\| {\bm e} \right\|}}\\
~~~ =  - \frac{1}{{{T_c}}}\frac{{{{\bm e}^T}}}{{\left\| {\bm e} \right\|}}{\left[ {\frac{{\partial W({{\left\| {\bm e} \right\|}^{m_0}})}}{{\partial {\bm e}}}\frac{{{m_0}{{\left\| {\bm e} \right\|}^{{m_0} - 1}}{{\bm e}^T}}}{{\left\| {\bm e} \right\|}}} \right]^{ - 1}}\\
~~~ =  - \frac{1}{{{m_0}{T_c}}}{\left[ {\frac{{\partial W({{\left\| {\bm e} \right\|}^{m_0}})}}{{\partial {{\left\| {\bm e} \right\|}^{m_0}}}}} \right]^{ - 1}}{\left\| {\bm e} \right\|^{1 - {m_0}}}\\
~~~ =  - \frac{1}{{{m_0}{T_c}}}{\left[ {\frac{{\partial W({V_2}^{m_0})}}{{\partial {V_2}^{m_0}}}} \right]^{ - 1}}{V_2}^{1 - {m_0}}
\end{array}\label{eq:107}\end{equation}

Then the reaching time of $\bm e = \bm 0$ can be deduced as:
\begin{equation}\begin{array}{l}
T({\bm x_0}) = \int\limits_{{V_2}({\bm x_0})}^{{V_2}(\bm 0)} {\frac{1}{{ - \frac{1}{{{m_0}{T_c}}}{{\left[ {\frac{{\partial W({V_2}^{m_0})}}{{\partial {V_2}^{m_0}}}} \right]}^{ - 1}}{V_2}^{1 - {m_0}}}}d} {V_2}\\
 ~~~~~~~~~= {T_c}\int\limits_{{V_2}(\bm 0)}^{{V_2}({\bm x_0})} {\frac{{\partial W({V_2}^{m_0})}}{{\partial {V_2}^{m_0}}}d} {V_2}^{m_0}\\
 ~~~~~~~~~= {T_c}(W({V_2}{({\bm x_0})^{m_0}}) - W({V_2}{(\bm 0)^{m_0}}))
\end{array}\label{eq:108}\end{equation}

Then the minimum upper bound of settling time of system can be predefined as $n{T_c}$.

The system will be said to be predefined-time stable if $n$ is finite.

The proof is complete.

According to general Leibniz rule, \eqref{eq:40} can be written as:
\begin{equation}\begin{array}{l}
\bm u =  - {b^{ - 1}}(\bm x,t)\left[ {f(\bm x,t) + \frac{1}{{{m_k}{T_c}}}\sum\limits_{k = 0}^{n - 2} {\sum\limits_{j = 0}^{n - k - 1} {\left( {\begin{array}{*{20}{c}}
{n - k - 1}\\
j
\end{array}} \right) }} } \right.\\
~~~~~~~\cdot {{({{\left[ {\frac{{\partial W({{\left\| {{{\bm s}_k}} \right\|}^{m_k}})}}{{\partial {{\left\| {{{\bm s}_k}} \right\|}^{m_k}}}}} \right]}^{ - 1}})}^{(n - k - 1 - j)}}{{(\frac{{{{\bm s}_k}}}{{{{\left\| {{{\bm s}_k}} \right\|}^{m_k}}}})}^{(j)}} \\
~~~~~~\left. { + \frac{1}{{{m_{n-1}}{T_c}}}{{\left[ {\frac{{\partial W({{\left\| {{{\bm s}_{n - 1}}} \right\|}^{m_{n-1}}})}}{{\partial {{\left\| {{{\bm s}_{n - 1}}} \right\|}^{m_{n-1}}}}}} \right]}^{ - 1}}\frac{{{{\bm s}_{n - 1}}}}{{{{\left\| {{{\bm s}_{n - 1}}} \right\|}^{m_{n-1}}}}}} \right]\\
~~~ =  - {b^{ - 1}}(\bm x,t)\left\{ {f(\bm x,t) + \frac{1}{{{m_k}{T_c}}}\sum\limits_{k = 0}^{n - 2} {\sum\limits_{j = 0}^{n - k - 1} {\left( {\begin{array}{*{20}{c}}
{n - k - 1}\\
j
\end{array}} \right)} } } \right.\\
\cdot {{({{\left[ {\frac{{\partial W({{\left\| {{{\bm s}_k}} \right\|}^{m_k}})}}{{\partial {{\left\| {{{\bm s}_k}} \right\|}^{m_k}}}}} \right]}^{ - 1}})}^{(n - k - 1 - j)}}(\sum\limits_{l = 0}^j {\left( {\begin{array}{*{20}{c}}
j\\
l
\end{array}} \right){{({\bm s_k})}^{(j - l)}}{{({{\left\| {{{\bm s}_k}} \right\|}^{ - {m_k}}})}^{(l)}}} ) \\
~~~~~~~\left. { + \frac{1}{{{m_{n-1}}{T_c}}}{{\left[ {\frac{{\partial W({{\left\| {{{\bm s}_{n - 1}}} \right\|}^{m_{n-1}}})}}{{\partial {{\left\| {{{\bm s}_{n - 1}}} \right\|}^{m_{n-1}}}}}} \right]}^{ - 1}}\frac{{{{\bm s}_{n - 1}}}}{{{{\left\| {{{\bm s}_{n - 1}}} \right\|}^{m_{n-1}}}}}} \right\}

\end{array}\label{eq:41}\end{equation}

\textbf{Result 2:} According to Result 1, the \textbf{singularity}, which occurs at the reaching time of sliding mode $\bm s_k$: $t\in ((n-k-1)T_c,(n-k)T_c)$, can also be \textbf{avoided} if $m_k=1$ for any $k$ and Assumption 2 is satisfied.

Proof of Result 2 (for nonsingular parameter selection): In this proof, we will study how to make control input finite when $\bm e \to \bm 0$ and $\left\| \bm \Im \right\| \ne 0$

Consider $m_k=1$ for all $k$, then a term in the first equation of \eqref{eq:41} can be expressed:
\begin{equation}\mathop {\lim }\limits_{{\bm s_k} \to 0} \frac{{{\bm s_k}}}{{{{\left\| {{\bm s_k}} \right\|}}}} = \frac{\bm 0}{{{{\left\| \bm 0 \right\|}}}} = {\bm E_{n1}}.\label{eq:109}\end{equation}
If ${\bm s_k} \to \bm 0$, the control input $\bm u$ can be written as:
\begin{equation}\begin{array}{l}
\bm u =  - {b^{ - 1}}(\bm x,t)\left[ {f(\bm x,t)+ \frac{1}{{{m_k}{T_c}}}\sum\limits_{k = 0}^{n - 2} {\sum\limits_{j = 0}^{n - k - 1} {\left( {\begin{array}{*{20}{c}}
{n - k - 1}\\
j
\end{array}} \right)} } } \right.\\
~~~~~~~\cdot {{({{\left[ {\frac{{\partial W({{\left\| {{{\bm s}_k}} \right\|}})}}{{\partial {{\left\| {{{\bm s}_k}} \right\|}}}}} \right]}^{ - 1}})}^{(n - k - 1 - j)}}{{({\bm E_{n1}})^{(j)}}} \\
~~~~~~\left. { + \frac{1}{{{T_c}}}{{\left[ {\frac{{\partial W({{\left\| {{{\bm s}_{n - 1}}} \right\|}^{}})}}{{\partial {{\left\| {{{\bm s}_{n - 1}}} \right\|}^{}}}}} \right]}^{ - 1}} \bm E_{n1}} \right]\\
\end{array}\label{eq:110}\end{equation}

It is clear that $\bm E_{n1}^{(j)}$ will not be infinite forever. Then, by considering a Faa di Bruno's formula in \cite{encinas2003short}, the singularity problem in controller \eqref{eq:41} will be focused on the following:
\begin{equation}\begin{array}{l}
{({\left[ {\frac{{\partial W({{\left\| {{{\bm s}_k}} \right\|}})}}{{\partial {{\left\| {{{\bm s}_k}} \right\|}}}}} \right]^{ - 1}})^{(n - k - 1 - j)}}\\
 = \sum {\frac{{(n - k - 1 - j)!}}{{{b_1}!{b_2}! \cdots {b_{n - k - 1 - j}}!}}\frac{{{\partial ^{{h_p}}}{{\left[ {\frac{{\partial W({{\left\| {{{\bm s}_k}} \right\|}})}}{{\partial {{\left\| {{{\bm s}_k}} \right\|}}}}} \right]}^{ - 1}}}}{{\partial {{({{\left\| {{{\bm s}_k}} \right\|^2}})}^{{h_p}}}}}}\\
 ~~~~\cdot {{(\frac{{{{({{\left\| {{{\bm s}_k}} \right\|^2}})}^{(1)}}}}{{1!}})}^{{b_1}}}{{(\frac{{{{(\left\| {{{\bm s}_k}} \right\|^2)}^{(2)}}}}{{2!}})}^{{b_2}}} \cdots {{(\frac{{{{(\left\| {{{\bm s}_k}} \right\|^2)}^{(n - k - 1 - j)}}}}{{(n - k - 1 - j)!}})}^{{b_{n - k - 1 - j}}}}
\end{array}\label{eq:111}\end{equation}
where $b_1+2b_2+\cdots+(n - k - 1 - j)b_{n - k - 1 - j}=n - k - 1 - j$ and $b_1+b_2+\cdots+b_{n - k - 1 - j}=h_p$ in which $b$ is non-negative integer.

\textbf{Remark 5:}; There is a new definition in Faa di Bruno's formula $\frac{{{\partial ^\gamma }\varphi (\phi (x))}}{{\partial {x^\gamma }}} = \sum {\frac{{\gamma !}}{{{b_1}! \cdots {b_\gamma }}}} \frac{{{\partial ^{{h_\gamma }}}\varphi (\phi (x))}}{{\partial \phi {{(x)}^{{h_\gamma }}}}}{(\frac{{\dot \phi (x)}}{{1!}})^{{b_1}}}$
$\cdots {(\frac{{{\phi ^{(\gamma )}}(x)}}{{\gamma !}})^{{b_\gamma }}}$ for this paper: $\frac{{{\partial ^\gamma }\varphi (\phi (x))}}{{\partial {x^\gamma }}}$ will be equal to ${\varphi (\phi (x))}$ if $\gamma$=0.

Let $\left\|{\bm {s_k}} \right\|^2={\bm s_k}^T{\bm s_k}$, then, $\left\|{\bm {s_k}} \right\|^2$ in \eqref{eq:111} can be expressed:
\begin{equation}\begin{array}{*{20}{l}}
{{{({\bm {s_k}^T}\bm {s_k})}^{({l_v})}} = \sum\limits_{v = 1}^n {{{(s_{k,v}^2)}^{({l_v})}}} }\\
~~~~~~~~~~~~~~~~{ = \sum\limits_{v = 1}^n {\sum {\left\{ {\frac{{{l_v}!}}{{{b_1}!{b_2}! \cdots {b_{{l_v}}}!}}\left[ {\frac{{{\partial ^{{h_v}}}s_{k,v}^2}}{{\partial {{({s_{k,v}})}^{{h_v}}}}}} \right]} \right.} } }\\
~~~~~~~~~~~~~~~~~~~{ \cdot \left. {{{(\frac{{{{\dot s}_{k,v}}}}{{1!}})}^{{b_1}}}{{(\frac{{{{\ddot s}_{k,v}}}}{{2!}})}^{{b_2}}} \cdots {{(\frac{{{{({s_{k,v}})}^{(l_v)}}}}{{v!}})}^{{b_{{l_v}}}}}} \right\}}
\end{array}\label{eq:114}\end{equation}
where ${\bm s_k} = {\left[ {\begin{array}{*{20}{c}}
{{s_{k,1}}}&{{s_{k,2}}}& \cdots &{{s_{k,n}}}
\end{array}} \right]^T}$; ${\bm e} = {\left[ {\begin{array}{*{20}{c}}
{{e_{1}}}&{{e_{2}}}& \cdots &{{e_{n}}}
\end{array}} \right]^T}$; $l_v=1,2\cdots,l_1$; $b_1+2b_2+\cdots+{l_v}b_{l_v}={l_v}$ and $b_1+b_2+\cdots+b_{l_v}=h_v$.

One of the possible singularity term in \eqref{eq:114} can be focused on the following:
\begin{equation}{\left[ {\frac{{{\partial ^{{h_v}}}s_{k,v}^2}}{{\partial {{({s_{k,v}})}^{{h_v}}}}}} \right] = \left\{ {\begin{array}{*{20}{l}}
{\sum\limits_{v = 1}^n {\left|2{s_{k,v}}\right|,{h_v} = 1} }\\
{\sum\limits_{v = 1}^n {2~~~~~~~~,{h_v} = 2} }\\
{~0~~~~~~~~~~~~,{h_v} \ge 3}
\end{array}} \right.}\label{eq:115}\end{equation}
and it is clear that the term \eqref{eq:115} will not be infinite if $s_{k,v}$ is not infinite.

\textbf{Assumption 3:} There will be no infinite value in \eqref{eq:115} if $s_{k,v}$ is bounded.

Then, other possible singularity terms: $\bm s_{k}^{(i_k)}$ for $i_k = 1, 2 \cdots, n-k-1-j$ in \eqref{eq:114} will be discussed in the following steps:

Step 0: \textbf{The singularity definition of this subsection is that control input will be infinite if $\bm e \to \bm 0$ and $\left\| \bm \Im \right\| \ne 0$}. $\bm s_{1,2,\cdots,n-1}^{(i_1)}$ will be assumed to be bounded first for $i_1 = 1, 2, \cdots, n-j$. Then, by substituting $\bm s_0=\bm e$ in \eqref{eq:111} and \eqref{eq:114}, it is clear that the singularity of controller about $\bm s_0$ will not occur if $m_k=1$ for any $k$.

Step 1: $\bm s_{2,3,\cdots,n-1}^{(i_1)}$ will be assumed to be bounded first for $i_1 = 1, 2, \cdots, n-1-j$. If $m_0=1$ and $\bm s_0=\bm e$, the singularity of controller about $\bm s_1$ can be expressed as:
\begin{equation}\begin{array}{l}
\begin{array}{*{20}{l}}
{\bm s_1^{({i_1})} = {\bm e^{({i_1})}} + \frac{1}{{{}{T_c}}}\sum\limits_{{j_1} = 0}^{{i_1}} {\left( {\begin{array}{*{20}{c}}
{{i_1}}\\
{{j_1}}
\end{array}} \right)} }\\
{ \cdot {{({{\left[ {\frac{{\partial W(\left\| \bm e \right\|)}}{{\partial \left\|\bm  e \right\|}}} \right]}^{ - 1}})}^{({i_1} - {j_1})}}{{({E_{n1}})}^{({j_1})}}}
\end{array}\\
 = {\bm e^{({i_1})}} + \frac{1}{{{}{T_c}}}\sum\limits_{{j_1} = 0}^{{i_1}} {\left( {\begin{array}{*{20}{c}}
{{i_1}}\\
{{j_1}}
\end{array}} \right)} \\
\left[ {\sum {\frac{{({i_1} - {j_1})!}}{{{b_1}!{b_2}! \cdots {b_{{i_1} - {j_1}}}!}}\frac{{{\partial ^{{h_1}}}{{\left[ {\frac{{\partial W(\left\| \bm e \right\|)}}{{\partial \left\| \bm e \right\|}}} \right]}^{ - 1}}}}{{\partial {{(\left\| \bm e \right\|)}^{{h_1}}}}}} } \right.{(\frac{{{{(\left\| \bm e \right\|)}^{(1)}}}}{{1!}})^{{b_1}}}\\
 \cdot \left. {{{(\frac{{{{(\left\| \bm e \right\|)}^{(2)}}}}{{2!}})}^{{b_2}}} \cdots {{(\frac{{{{(\left\| \bm e \right\|)}^{({i_1} - {j_1})}}}}{{({i_1} - {j_1})!}})}^{{b_{{i_1} - {j_1}}}}}} \right]{({\bm E_{n1}})^{({j_1})}}
\end{array}\label{eq:116}\end{equation}
where $b_1+2b_2+\cdots+(i_1-j_1)b_{i_1-j_1}=i_1-j_1$ and $b_1+b_2+\cdots+b_{i_1-j_1}=h_{1}$.

Let $\left\|{\bm e} \right\|^2={\bm e}^T{\bm e}$, then:
\begin{equation}\begin{array}{*{20}{l}}
{{{({\bm e^T}\bm e)}^{({l_{v1}})}} = \sum\limits_{{v1} = 1}^n {{{(e_{v1}^2)}^{({l_{v1}})}}} }\\
{ = \sum\limits_{{v1} = 1}^n {\sum {\left\{ {\frac{{{l_{v1}}!}}{{{b_1}!{b_2}! \cdots {b_{{l_{v1}}}}!}}\left[ {\frac{{{\partial ^{{h_{v1}}}}e_{v1}^2}}{{\partial {{({e_{v1}})}^{{h_{v1}}}}}}} \right]} \right.} } }\\
{ \cdot \left. {{{(\frac{{{{\dot e}_{v1}}}}{{1!}})}^{{b_1}}}{{(\frac{{{{\ddot e}_{v1}}}}{{2!}})}^{{b_2}}} \cdots {{(\frac{{{{({e_{v1}})}^{({v1})}}}}{{{v1}!}})}^{{b_{{l_{v1}}}}}}} \right\}}
\end{array}\label{eq:117}\end{equation}
where ${\bm e} = {\left[ {\begin{array}{*{20}{c}}
{{e_{v1}}}&{{e_{v2}}}& \cdots &{{e_{vn}}}
\end{array}} \right]^T}$; $l_{v1}=1,2\cdots,l_1$; $b_1+2b_2+\cdots+{l_{v1}}b_{l_{v1}}={l_{v1}}$ and $b_1+b_2+\cdots+b_{l_{v1}}=h_{v1}$.

The singular problem will be focused on the following:
\begin{equation}{\left[ {\frac{{{\partial ^{{h_{v1}}}}e_{v1}^2}}{{\partial {{({e_{v1}})}^{{h_{v1}}}}}}} \right] = \left\{ {\begin{array}{*{20}{l}}
{\sum\limits_{{v1} = 1}^n {\left|2{e_{v1}}\right|,{h_{v1}} = 1} }\\
{\sum\limits_{{v1} = 1}^n {2,{h_{v1}} = 1} }\\
{0,{h_{v1}} \ge 3}
\end{array}} \right.}\label{eq:118}\end{equation}

Then, by considering \eqref{eq:116}, \eqref{eq:117} and \eqref{eq:118}, The singularity of controller about $s_1$ will not occur if $m_0=1$.

Steps $2, \cdots ,k-1$: ~~~~~~~~~~~~~~~$\vdots $

Step k: $\bm s_{k+1,\cdots,n-1}^{(i_2)}$ will be be assumed to be bounded first for $i_k = 1,2, \cdots ,n-k-j$. The singularity of controller about $\bm s_k$ if $m_{k-1}=1$  can be expressed as:
\begin{equation}\begin{array}{*{20}{l}}
{\begin{array}{*{20}{l}}
{\bm s_k^{({i_k})} = {\bm s_{k-1}}^{({i_k})} + \frac{1}{{{T_c}}}\sum\limits_{{j_k} = 0}^{{i_k}} {\left( {\begin{array}{*{20}{c}}
{{i_k}}\\
{{j_k}}
\end{array}} \right)} }\\
{ \cdot {{({{\left[ {\frac{{\partial W(\left\| {{\bm s_{k-1}}} \right\|)}}{{\partial \left\| {{\bm s_{k-1}}} \right\|}}} \right]}^{ - 1}})}^{({i_k} - {j_k})}}{{({\bm E_{n1}})}^{({j_k})}}}
\end{array}}\\
{ = {\bm s_{k-1}}^{({i_k})} + \frac{1}{{{T_c}}}\sum\limits_{{j_k} = 0}^{{i_k}} {\left( {\begin{array}{*{20}{c}}
{{i_k}}\\
{{j_k}}
\end{array}} \right)} }\\
{\left[ {\sum {\frac{{({i_k} - {j_k})!}}{{{b_1}!{b_2}! \cdots {b_{{i_k} - {j_k}_{}}}!}}\frac{{{\partial ^{{h_k}}}{{\left[ {\frac{{\partial W(\left\| {{\bm s_{k-1}}} \right\|)}}{{\partial \left\| {{\bm s_{k-1}}} \right\|}}} \right]}^{ - 1}}}}{{\partial {{(\left\| {{\bm s_{k-1}}} \right\|)}^{{h_k}}}}}} } \right.{{(\frac{{{{(\left\| {{\bm s_{k-1}}} \right\|)}^{(1)}}}}{{1!}})}^{{b_1}}}}\\
{ \cdot \left. {{{(\frac{{{{(\left\| {{\bm s_{k-1}}} \right\|)}^{(2)}}}}{{2!}})}^{{b_2}}} \cdots {{(\frac{{{{(\left\| {{\bm s_{k-1}}} \right\|)}^{({i_k} - {j_k})}}}}{{({i_k} - {j_k})!}})}^{{b_{{i_k} - {j_k}}}}}} \right]{{({\bm E_{n1}})}^{({j_k})}}}
\end{array}\label{eq:119}\end{equation}
where $b_1+2b_2+\cdots+(i_k-j_k)b_{i_k-j_k}=i_k-j_k$ and $b_1+b_2+\cdots+b_{i_k-j_k}=h_{k}$.

Let $\left\|{\bm {s_{k-1}}} \right\|^2={\bm {s_{k-1}}}^T{\bm {s_{k-1}}}$, then:
\begin{equation}\begin{array}{*{20}{l}}
{{{({\bm {s_{k-1}}^T}\bm {s_{k-1}})}^{({l_{vk}})}} = \sum\limits_{{vk}= 1}^n {{{({s_{{k-1},{vk}}}^2)}^{({l_{vk}})}}} }\\
{ = \sum\limits_{{vk} = 1}^n {\sum {\left\{ {\frac{{{l_{vk}}!}}{{{b_1}!{b_2}! \cdots {b_{{l_{vk}}}}!}}\left[ {\frac{{{\partial ^{{h_{vk}}}}s_{{k-1},{vk}}^2}}{{\partial {{({{s_{{k-1},vk}}})}^{{h_{v2}}}}}}} \right]} \right.} } }\\
{ \cdot \left. {{{(\frac{{{{ {\dot s_{{k-1},{vk}}}}}}}{{1!}})}^{{b_1}}}{{(\frac{{{{ {\ddot s_{{k-1},{vk}}}}}}}{{2!}})}^{{b_2}}} \cdots {{(\frac{{{{({{s_{{k-1},{vk}}}})}^{({vk})}}}}{{{vk}!}})}^{{b_{{l_{vk}}}}}}} \right\}}
\end{array}\label{eq:120}\end{equation}
where ${\bm s_{k-1}} = {\left[ {\begin{array}{*{20}{c}}
{{s_{{k-1},v1}}}&{{s_{{k-1},v2}}}& \cdots &{{s_{{k-1},vn}}}
\end{array}} \right]^T}$; $l_{vk}=1,2\cdots,l_k$; $b_1+2b_2+\cdots+{l_{vk}}b_{l_{vk}}={l_{vk}}$ and $b_1+b_2+\cdots+b_{l_{vk}}=h_{vk}$.

with
\begin{equation}{\left[ {\frac{{{\partial ^{{h_{vk}}}}s_{{k-1},vk}^2}}{{\partial {{({s_{{k-1},vk}})}^{{h_{vk}}}}}}} \right] = \left\{ {\begin{array}{*{20}{l}}
{\sum\limits_{{vk} = 1}^n {\left|{2s_{{k-1},vk}}\right|,{h_{vk}} = 1} }\\
{\sum\limits_{{vk} = 1}^n {2,{h_{vk}} = 2} }\\
{0,{h_{vk}} \ge 3}
\end{array}} \right.}\label{eq:121}\end{equation}

According to step 1 (\eqref{eq:118} is substituted into \eqref{eq:117}, and \eqref{eq:117} is substituted into \eqref{eq:116}), it is clear that $s_{1,vk}^{(1),(2)\cdots,(vk)}$ will not be infinite if $m_{0,1,\cdots k-2}=1$. Then, similar to step 1, by considering \eqref{eq:119}, \eqref{eq:120} and \eqref{eq:121}, The singularity of controller about $\bm s_k$ will also not occur if $m_{k-1}= 1$.

Steps $k+1 \cdots ,n-1$: ~~~~~~~~~~~~~~~$\vdots $

Step $n$: Make an identity: $\bm s_n =\bm 0$ in \eqref{eq:38}. According to the above-mentioned $n-1$ steps, $\bm s_{1,2,\cdots,n-2}^{(i_2)}$ can be predefined to be bounded if $m_{0, 1 \cdots ,n-2}=1$. The singularity of controller about $\bm s_{n}=\bm 0$ if $m_{n-1}=1$ can be expressed as:
\begin{equation}\begin{array}{l}
\frac{1}{{{m_{n - 1}}{T_c}}}{\left[ {\frac{{\partial W({{\left\| {{\bm s_{n - 1}}} \right\|}^{{m_{n - 1}}}})}}{{\partial {{\left\| {{\bm s_{n - 1}}} \right\|}^{{m_{n - 1}}}}}}} \right]^{ - 1}}\frac{{{\bm s_{n - 1}}}}{{{{\left\| {{\bm s_{n - 1}}} \right\|}^{{m_{n - 1}}}}}}\\
 = \frac{1}{{{T_c}}}{\left[ {\frac{{\partial W(\left\| {{\bm s_{n - 1}}} \right\|)}}{{\partial \left\| {{\bm s_{n - 1}}} \right\|}}} \right]^{ - 1}}{\bm E_{n1}}
\end{array}\label{eq:122}\end{equation}

It is clear that the singularity of controller about $s_{n}$ will also not occur if $m_{n-1}= 1$.

\textbf{According to the above steps, the Assumption 3 is actually satisfied everywhere if $m_k=1$}. Hence, the control input will not infinite if $m_k=1$.

The proof is complete.

However, according to the simulation in Fig.3 and Result 1. It is clear that the predefined-time stability cannot be guaranteed if $m_k=1$ for everywhere. Hence, a switching decision similar to \eqref{eq:14} will be discussed in the following.

\textbf{Result 3:}
In order to ensure the predefined-time stability of the system, parameter switching can only occur on the sliding mode surface $\bm s_{k+1}$ relative to virtual system state $\bm s_{k}$. Hence, in reaching stage of virtual system state $\bm s_{k}$, the corresponding parameter switching cannot be made. \textbf{In addition, the selection range of nonsingular parameters should be: $0<m_k\le \frac{1}{n-k}$} and $m_k=1$

Proof of the selection range of nonsingular parameters Result 3 ($0<m_k\le \frac{1}{n-k}$):

According to the Faa di Bruno's formula, \eqref{eq:41} can be written as:
\begin{equation}\begin{array}{*{20}{l}}
{\bm u =  - {b^{ - 1}}(\bm x,t)\left\{ {f(\bm x,t) + \frac{1}{{{m_k}{T_c}}}\sum\limits_{k = 0}^{n - 2} {\sum\limits_{j = 0}^{n - k - 1} {\left( {\begin{array}{*{20}{c}}
{n - k - 1}\\
j
\end{array}} \right)} } } \right.}\\
{ \cdot \left[ {(\sum \frac{{(n - k - 1 - j)!}}{{{b_1}!{b_2}! \cdots {b_{n - k - 1 - j}}!}}\frac{{{\partial ^{{h_p}}}{{\left[ {\frac{{\partial W({{\left\| {{\bm s_k}} \right\|}^{{m_k}}})}}{{\partial {{\left\| {{\bm s_k}} \right\|}^{{m_k}}}}}} \right]}^{ - 1}}}}{{\partial {{({{\left\| {{\bm s_k}} \right\|}^{{m_k}}})}^{{h_p}}}}}}{(\frac{{{{({{\left\| {{\bm s_k}} \right\|}^{{m_k}}})}^{(1)}}}}{{1!}})^{{b_1}}} \right.}\\
\begin{array}{l}
{(\frac{{{{({{\left\| {{\bm s_k}} \right\|}^{{m_k}}})}^{(2)}}}}{{2!}})^{{b_2}}} \cdots {(\frac{{{{({{\left\| {{\bm s_k}} \right\|}^{{m_k}}})}^{(n - k - 1 - j)}}}}{{n - k - 1 - j!}})^{{b_{n - k - 1 - j}}}})(\sum\limits_{l = 0}^j \left( {\begin{array}{*{20}{c}}
j\\
l
\end{array}} \right)\\
 \cdot {\left[ {{{({\bm s_k})}^{(j - l)}}\sum (\frac{{l!}}{{{b_1}!{b_2}! \cdots {b_l}!}}(\prod\limits_{z = 0}^{h - 1} {( - \frac{{{m_k}}}{2} - z))} }{{\left\| {{\bm s_k}} \right\|}^{ - {m_k} - 2h}} \right.}
\end{array}\\
{ \cdot \left. {(\left. {\frac{{{{({{\left\| {{\bm s_k}} \right\|}^2})}^{(1)}}}}{{1!}}{)}^{{b_1}}{{(\frac{{{{({{\left\| {{\bm s_k}} \right\|}^2})}^{(2)}}}}{{2!}})}^{{b_2}}} \cdots {{(\frac{{{{({{\left\| {{\bm s_k}} \right\|}^2})}^{(l)}}}}{{l!}})}^{{b_l}}})} \right])} \right]}\\
{\left. { + \frac{1}{{{m_{n - 1}}{T_c}}}{{\left[ {\frac{{\partial W({{\left\| {{\bm s_{n - 1}}} \right\|}^{{m_{n - 1}}}})}}{{\partial {{\left\| {{\bm s_{n - 1}}} \right\|}^{{m_{n - 1}}}}}}} \right]}^{ - 1}}\frac{{{\bm s_{n - 1}}}}{{{{\left\| {{\bm s_{n - 1}}} \right\|}^{{m_{n - 1}}}}}}} \right\}}
\end{array}\label{eq:42}\end{equation}
with
\begin{equation}\begin{array}{l}
{({\left\| {{\bm s_k}} \right\|^{{m_k}}})^{^{(n - k - 1 - j)}}}\\
 = \sum \frac{{(n - k - 1 - j)!}}{{{b_1}!{b_2}! \cdots {b_{n - k - 1 - j}}!}}\left[ {\prod\limits_{z = 0}^{{h_p} - 1} {(  \frac{{{m_k}}}{2} - z)} } \right]{\left\| {{\bm s_k}} \right\|^{ {m_k} - 2{h_p}}}\\
 \cdot {(\frac{{{{({{\left\| {{\bm s_k}} \right\|}^2})}^{(1)}}}}{{1!}})^{{b_1}}}{(\frac{{{{({{\left\| {{\bm s_k}} \right\|}^2})}^{(2)}}}}{{2!}})^{{b_2}}} \cdots {(\frac{{{{({{\left\| {{\bm s_k}} \right\|}^2})}^{(n - k - 1 - j)}}}}{{(n - k - 1 - j)!}})^{{b_{n - k - 1 - j}}}}
\end{array}\label{eq:123}\end{equation}
where $b_1+2b_2+\cdots+lb_l=l$ and $b_1+b_2+\cdots+b_l=h$; $b_1+2b_2+\cdots+(n - k - 1 - j)b_{n - k - 1 - j}=n - k - 1 - j$ and $b_1+b_2+\cdots+b_{n - k - 1 - j}=h_p$.

Let $\left\|{\bm s}_k \right\|^2={\bm s}_k^T{\bm s}_k$, then:
\begin{equation}\begin{array}{l}
{({\bm s}_k^T{{\bm s}_k})^{({l_v})}} = \sum\limits_{v = 1}^n {{{({ s}_{k,v}^2)}^{({l_v})}}} \\
~~~~~~~~~~~~~~~ =\sum\limits_{v = 1}^n {\sum {\left\{ {\frac{{{l_v}!}}{{{b_1}!{b_2}! \cdots {b_{{l_v}}}!}}\left[ {\frac{{{\partial ^{{h_v}}}s_{k,v}^2}}{{\partial {{({s_{k,v}})}^{{h_v}}}}}} \right]} \right.} } \\
~~~~~~~~~~~~~~~~~~~~ \cdot \left. {{{(\frac{{{{\dot s}_{k,v}}}}{{1!}})}^{{b_1}}}{{(\frac{{{{\ddot s}_{k,v}}}}{{2!}})}^{{b_2}}} \cdots {{(\frac{{{{({s_{k,v}})}^{(v)}}}}{{v!}})}^{{b_{{l_v}}}}}} \right\}
\end{array}\label{eq:43}\end{equation}
where ${\bm s_k} = {\left[ {\begin{array}{*{20}{c}}
{{s_{k,1}}}&{{s_{k,2}}}& \cdots &{{s_{k,n}}}
\end{array}} \right]^T}$. $b_1+2b_2+\cdots+{l_v}b_{l_v}={l_v}$ and $b_1+b_2+\cdots+b_{l_v}=h_v$.

Step 0: By assuming that $\bm s_0 \ne \bm 0$ and $\bm s_{1,2 \cdots, n-1}=\bm 0$, then, according to the results of Lemma1 and \eqref{eq:38}, $( {{{\dot {\bm s}}_{1}},{{ \ddot {\bm s}}_{1}}, \cdots ,{{({{\bm s}_{1}})}^{(n-1)}}} )$ have to converge to zero in $t \ge (n-1)T_c$ if the system is predefined-time stable and $0<m_0\le \frac{1}{n}$. Then, the control input $\bm u$ can be written as:
\begin{equation}\begin{array}{*{20}{l}}
{\bm u =  - {b^{ - 1}}(\bm x,t)\left\{ {f(\bm x,t) + \frac{1}{{{m_0}{T_c}}}\sum\limits_{j = 0}^{n - 1} {\left( {\begin{array}{*{20}{c}}
{n - 1}\\
j
\end{array}} \right)} } \right.}\\
{ \cdot \left[ {(\sum \frac{{(n - 1 - j)!}}{{{b_1}!{b_2}! \cdots {b_{n - 1 - j}}!}}\frac{{{\partial ^{{h_p}}}{{\left[ {\frac{{\partial W({{\left\|\bm  e \right\|}^{{m_0}}})}}{{\partial {{\left\| \bm e \right\|}^{{m_0}}}}}} \right]}^{ - 1}}}}{{\partial {{({{\left\| \bm e \right\|}^{{m_0}}})}^{{h_p}}}}}} \right.}\\
\begin{array}{l}
{(\frac{{{{({{\left\| \bm e \right\|}^{{m_0}}})}^{(1)}}}}{{1!}})^{{b_1}}}{(\frac{{{{({{\left\| \bm e \right\|}^{{m_0}}})}^{(2)}}}}{{2!}})^{{b_2}}} \cdots {(\frac{{{{({{\left\| \bm e \right\|}^{{m_0}}})}^{(n - 1 - j)}}}}{{(n - 1 - j)!}})^{{b_{n - 1 - j}}}})\\
\cdot (\sum\limits_{l = 0}^j {\left( {\begin{array}{*{20}{c}}
j\\
l
\end{array}} \right)\left[ {{{(\bm e)}^{(j - l)}}\sum (\frac{{l!}}{{{b_1}!{b_2}! \cdots {b_l}!}}(\prod\limits_{z = 0}^{h - 1} {( - \frac{{{m_0}}}{2} - z))} } \right.}
\end{array}\\
{ \cdot {{\left\| \bm e \right\|}^{ - {m_0} - 2h}}\left. {(\left. {\frac{{{{({{\left\| \bm e \right\|}^2})}^{(1)}}}}{{1!}}{)^{{b_1}}}{{(\frac{{{{({{\left\| \bm e \right\|}^2})}^{(2)}}}}{{2!}})}^{{b_2}}} \cdots {{(\frac{{{{({{\left\| \bm e \right\|}^2})}^{(l)}}}}{{l!}})}^{{b_l}}})} \right])} \right]}\\
{\left. { + \frac{1}{{{m_{n - 1}}{T_c}}}{{\left[ {\frac{{\partial W({{\left\| {{\bm s_{n - 1}}} \right\|}^{{m_{n - 1}}}})}}{{\partial {{\left\| {{\bm s_{n - 1}}} \right\|}^{{m_{n - 1}}}}}}} \right]}^{ - 1}}\frac{{{\bm s_{n - 1}}}}{{{{\left\| {{\bm s_{n - 1}}} \right\|}^{{m_{n - 1}}}}}}} \right\}}
\end{array}\label{eq:124}\end{equation}

According to \eqref{eq:38}, $\dot{\bm s_0} $ can be replaced by a function about $\bm s_0$ when $\bm s_1=\bm 0$:
\begin{equation}\begin{array}{l}
\bm {\dot s_0} =  - \frac{1}{{{T_c}}}{\left\| {\frac{{\partial W({{\left\| {{\bm s_0}} \right\|}^{{m_0}}})}}{{\partial {\bm s_0}}}} \right\|^{ - 2}}{\left[ {\frac{{\partial W({{\left\| {{\bm s_0}} \right\|}^{{m_0}}})}}{{\partial {\bm s_0}}}} \right]^T}\\
~~~ = \frac{1}{{{m_0}{T_c}}}{\left[ {\frac{{\partial W({{\left\| \bm e \right\|}^{{m_0}}})}}{{\partial {{\left\| \bm e \right\|}^{{m_0}}}}}} \right]^{ - 1}}\frac{\bm e}{{{{\left\| \bm e \right\|}^{{m_0}}}}}
\end{array}\label{eq:125}\end{equation}

Then, \eqref{eq:43} can be written by considering \eqref{eq:115} and \eqref{eq:125} as:
\begin{equation}\begin{array}{l}
\begin{array}{*{20}{l}}
{{{(\bm e^T{\bm e})}^{({l_v})}} = \sum\limits_{v = 1}^n {{{(e_v^2)}^{({l_v})}}} }\\
{ = \sum\limits_{v = 1}^n {\left\{ {\frac{{{l_v}!}}{{{b_1}!{b_2}! \cdots {b_{{l_v}}}!}}\left| {2{e_v}} \right|{{(\frac{1}{{{m_0}{T_c}}}{{\left[ {\frac{{\partial W({{\left\| {{e_v}} \right\|}^{{m_0}}})}}{{\partial {{\left\| {{e_v}} \right\|}^{{m_0}}}}}} \right]}^{ - 1}}e_v^{1 - {m_0}})}^{{b_1}}}} \right.} }\\
{{{\left. {{{(\frac{{{{\ddot e}_{v}}}}{{2!}})}^{{b_2}}} \cdots {{(\frac{{{{({e_v})}^{(v)}}}}{{v!}})}^{{b_{{l_v}}}}}} \right|}_{{h_v} = 1}} + 2\frac{{{l_v}!}}{{{b_1}!{b_2}! \cdots {b_{{l_v}}}!}}(\frac{1}{{{m_0}{T_c}}}}
\end{array}\\
\left. { \cdot {{\left. {{{\left[ {\frac{{\partial W({{\left\| {{e_v}} \right\|}^{{m_0}}})}}{{\partial {{\left\| {{e_v}} \right\|}^{{m_0}}}}}} \right]}^{ - 1}}e_v^{1 - {m_0}}{)^{{b_1}}}{{(\frac{{{{\ddot e}_{v}}}}{{2!}})}^{{b_2}}} \cdots {{(\frac{{{{({e_v})}^{(v)}}}}{{v!}})}^{{b_{{l_v}}}}}} \right|}_{{h_v} = 2}}} \right\}
\end{array}\label{eq:126}\end{equation}
with
\begin{equation}\left\{ {\begin{array}{*{20}{l}}
{{b_1} + {b_2} +  \cdots  + {b_{{l_v}}} = {h_v}}\\
{{b_1} + 2{b_2} +  \cdots  + {l_v}{b_{{l_v}}} = {l_v}}
\end{array}} \right.\label{eq:127}\end{equation}

Consider a singularity expression: $\bm e\to \bm 0$ with $\bm e^{(1,2 \cdots, n-1)}= a~bounded~vector$. According to \eqref{eq:115}, \eqref{eq:123} and \eqref{eq:126}, the main singular term (the highest-order negative exponential) of controller in \eqref{eq:124} about $\bm s_0=\bm e$ can be simplified by considering $\left\|\bm 0\right\|_1=\left\|\bm 0\right\|_2$ as:
\begin{equation}\begin{array}{*{20}{l}}
{\begin{array}{*{20}{l}}
{{{\left\| \bm e \right\|}^{\sum\limits_{v0 = 1}^{n - 1 - j} {{b_{v0}}\left[ {{m_0} - 2{h_p} + \sum\limits_{v0 = 1}^{n - 1 - j} {{b_{v0}}} (1 + {b_1}(1 - m_0))} \right]} }}{e^{(j-l)}}}\\
~~~~~~{ \cdot {{\left\| \bm e \right\|}^{ - {m_0} - 2h}}{{\left\| \bm e \right\|}^{\sum\limits_{v0 = 1}^l {{b_{v0}}(1 + {b_1}(1 - m_0))} }}}
\end{array}}\\
{\begin{array}{*{20}{l}}
{\mathop  = \limits^{\bm e \to 0} {{\left\| \bm e \right\|}^{{m_0}({h_p}(1 - {b_1}{h_p}) - 1 - h{b_1})}}}{  {{\left. {{{\left\| \bm e \right\|}^{(h_p^2 + h)({b_1} - 1) + 1}}{\bm E_{n1}}} \right|}_{{h_v} = 1,2}}}
\end{array}}
\end{array}\label{eq:128}\end{equation}
with
\begin{equation}{\left\{ {\begin{array}{*{20}{l}}
{\sum\limits_{v0 = 1}^{n - 1 - j} {{b_{v0}}} } = {h_p} \le {n-1-j}\\
~{\sum\limits_{v0 = 1}^{l} {{b_{v0}}} }~ = h \le {l}\\
~~{{h_v}- ({b_2} +  \cdots  + {b_{{l_v}})} = {b_1}} \le {l_v}
\end{array}} \right.}\label{eq:129}\end{equation}
where $j=l$ because the fianl exponential term of $\bm e^{(j-l)}$ will be smaller than $1$ if $j\ne l$ by considering the decomposition term $\sum\limits_{k = 0}^{n - 2} {(\frac{1}{{{T_c}}}{{\left\| {\frac{{\partial W({{\left\| {{\bm s_k}} \right\|}^{m_{k}}})}}{{\partial {\bm s_k}}}} \right\|}^{ - 2}}{{\left[ {\frac{{\partial W({{\left\| {{\bm s_k}} \right\|}^{m_{k}}})}}{{\partial {\bm s_k}}}} \right]}^T}} {)^{(n - k - 1)}}$ in \eqref{eq:40}.

Then, \eqref{eq:128} will not be infinite if the following inequality is satisfied:
\begin{equation}\begin{array}{l}
{m_0}( - {h_p}(1 - {b_1}{h_p}) + 1 + h{b_1}) \le (h_p^2 + h)({b_1} - 1) + (j - l)\\
 \Rightarrow \left\{ {\begin{array}{*{20}{c}}
~~{{m_0} \le (\frac{{(h_p^2 + h) + (j - l)}}{{({h_p}(2{h_p} - 1) + 1 + 2h)}})}\\
{{m_0} \le (\frac{{(j - l)}}{{({h_p}({h_p} - 1) + 1 + h)}})}
\end{array}} \right.\\
 \Rightarrow {m_0} \le (\frac{1}{{(1(1 - 1) + 1 + n - 1)}}) \le \frac{1}{n}
\end{array}\label{eq:130}\end{equation}
where $m_0$ will be minimum by considering the inequality relationship in \eqref{eq:129}.

Then the controller in condition of step 0 will not be singular if $0<m\le \frac{1}{n}$

Step $1,\cdots,k-1$:~~~~~~~~~~~~~~~~~\vdots

Step $k$: By assuming that $\bm s_{k+1,k+2 \cdots, n-1}=\bm 0$ and $\bm s_k \ne \bm 0$, then, according to the results of Lemma1 and \eqref{eq:38}, $( {{{\dot {\bm s}}_{k+1}},{{ \ddot {\bm s}}_{k+1}} \cdots {{({{\bm s}_{k+1}})}^{(n-k-1)}}} )$ have to converge to zero in $t\ge (n-k-1)T_c$ if the system is predefined-time stable and $0<m_k\le \frac{1}{n-k-1}$. Then, according to \eqref{eq:38}, it is clear that it is clear that . The control input $\bm u$ can be written as:
\begin{equation}\begin{array}{*{20}{l}}
{\bm u =  - {b^{ - 1}}(\bm x,t)\left\{ {f(\bm x,t) + \frac{1}{{{m_0}{T_c}}}\sum\limits_{j = 0}^{n - k - 1} {\left( {\begin{array}{*{20}{c}}
{n - k - 1}\\
j
\end{array}} \right)} } \right.}\\
{ \cdot \left[ {(\sum {\frac{{(n - k - 1 - j)!}}{{{b_1}!{b_2}! \cdots {b_{n - k - 1 - j}}!}}} \frac{{{\partial ^{{h_p}}}{{\left[ {\frac{{\partial W({{\left\| {{\bm s_k}} \right\|}^{{m_k}}})}}{{\partial {{\left\| {{\bm s_k}} \right\|}^{{m_k}}}}}} \right]}^{ - 1}}}}{{\partial {{({{\left\| {{\bm s_k}} \right\|}^{{m_k}}})}^{{h_p}}}}}} \right.}\\
{\begin{array}{*{20}{l}}
{{{(\frac{{{{({{\left\| \bm e \right\|}^{{m_k}}})}^{(1)}}}}{{1!}})}^{{b_1}}}{{(\frac{{{{({{\left\| \bm e \right\|}^{{m_k}}})}^{(2)}}}}{{2!}})}^{{b_2}}} \cdots {{(\frac{{{{({{\left\| \bm e \right\|}^{{m_k}}})}^{(n - 1 - j)}}}}{{(n - k-1 - j)!}})}^{{b_{n - k - 1 - j}}}})}\\
{ \cdot (\sum\limits_{l = 0}^j {\left( {\begin{array}{*{20}{c}}
j\\
l
\end{array}} \right)\left[ {{{(\bm e)}^{(j - l)}}\sum {(\frac{{l!}}{{{b_1}!{b_2}! \cdots {b_l}!}}(} \prod\limits_{z = 0}^{h - 1} {( - \frac{{{m_k}}}{2} - z))} } \right.} }
\end{array}}\\
{ \cdot {{\left\| \bm e \right\|}^{ - {m_k} - 2h}}\left. {(\left. {\frac{{{{({{\left\|\bm  e \right\|}^2})}^{(1)}}}}{{1!}}{)^{{b_1}}}{{(\frac{{{{({{\left\| \bm e \right\|}^2})}^{(2)}}}}{{2!}})}^{{b_2}}} \cdots {{(\frac{{{{({{\left\| \bm e \right\|}^2})}^{(l)}}}}{{l!}})}^{{b_l}}})} \right])} \right]}\\
{\left. { + \frac{1}{{{m_{n - 1}}{T_c}}}{{\left[ {\frac{{\partial W({{\left\| {{\bm s_{n - 1}}} \right\|}^{{m_{n - 1}}}})}}{{\partial {{\left\| {{\bm s_{n - 1}}} \right\|}^{{m_{n - 1}}}}}}} \right]}^{ - 1}}\frac{{{\bm s_{n - 1}}}}{{{{\left\| {{\bm s_{n - 1}}} \right\|}^{{m_{n - 1}}}}}}} \right\}}
\end{array}\label{eq:131}\end{equation}

According to \eqref{eq:38}, $\dot{\bm s_k} $ can be replaced by a function about $\bm s_0$ when $\bm s_{k+1}=\bm 0$:
\begin{equation}\begin{array}{l}
\bm {\dot s_k} =  - \frac{1}{{{T_c}}}{\left\| {\frac{{\partial W({{\left\| {{\bm s_k}} \right\|}^{{m_k}}})}}{{\partial {\bm s_k}}}} \right\|^{ - 2}}{\left[ {\frac{{\partial W({{\left\| {{\bm s_k}} \right\|}^{{m_k}}})}}{{\partial {\bm s_k}}}} \right]^T}\\
~~~ = \frac{1}{{{m_k}{T_c}}}{\left[ {\frac{{\partial W({{\left\| \bm s_k \right\|}^{{m_k}}})}}{{\partial {{\left\| \bm s_k \right\|}^{{m_k}}}}}} \right]^{ - 1}}\frac{\bm s_k}{{{{\left\| \bm s_k \right\|}^{{m_k}}}}}
\end{array}\label{eq:132}\end{equation}

Then, \eqref{eq:43} can be written by considering \eqref{eq:115} and \eqref{eq:132} as:
\begin{equation}\begin{array}{*{20}{l}}
{\begin{array}{*{20}{l}}
{{{(\bm s_k^T{\bm s_k})}^{({l_v})}} = \sum\limits_{v = 1}^n {{{(s_{k,v}^2)}^{({l_v})}}} }\\
{ = \sum\limits_{v = 1}^n {\left\{ {\frac{{{l_v}!}}{{{b_1}!{b_2}! \cdots {b_{{l_v}}}!}}\left| {2{s_{k,v}}} \right|{{(\frac{1}{{{m_k}{T_c}}}{{\left[ {\frac{{\partial W({{\left\| {{s_{k,v}}} \right\|}^{{m_k}}})}}{{\partial {{\left\| {{s_{k,v}}} \right\|}^{{m_k}}}}}} \right]}^{ - 1}}s_{k,v}^{1 - {m_k}})}^{{b_1}}}} \right.} }\\
{{{\left. {{{(\frac{{{{\ddot s}_{k,v}}}}{{2!}})}^{{b_2}}} \cdots {{(\frac{{{{({s_{k,v}})}^{(v)}}}}{{v!}})}^{{b_{{l_v}}}}}} \right|}_{{h_v} = 1}} + 2\frac{{{l_v}!}}{{{b_1}!{b_2}! \cdots {b_{{l_v}}}!}}(\frac{1}{{{m_k}{T_c}}}}
\end{array}}\\
{\left. { \cdot {{\left. {{{\left[ {\frac{{\partial W({{\left\| {{s_{k,v}}} \right\|}^{{m_k}}})}}{{\partial {{\left\| {{s_{k,v}}} \right\|}^{{m_k}}}}}} \right]}^{ - 1}}s_{k,v}^{1 - {m_k}}{)^{{b_1}}}{{(\frac{{{{\ddot s}_{k,v}}}}{{2!}})}^{{b_2}}} \cdots {{(\frac{{{{({s_{k,v}})}^{(v)}}}}{{v!}})}^{{b_{{l_v}}}}}} \right|}_{{h_v} = 2}}} \right\}}
\end{array}\label{eq:133}\end{equation}
with
\begin{equation}\left\{ {\begin{array}{*{20}{l}}
{{b_1} + {b_2} +  \cdots  + {b_{{l_v}}} = {h_v}}\\
{{b_1} + 2{b_2} +  \cdots  + {l_v}{b_{{l_v}}} = {l_v}}
\end{array}} \right.\label{eq:134}\end{equation}

Consider a singularity expression: $\bm e\to \bm 0$ with $\bm e^{(1,2 \cdots, n-1)}= a~bounded~constant$. According to \eqref{eq:115}, \eqref{eq:123} and \eqref{eq:133}, the main singular term (the highest-oder negative exponential) of controller in \eqref{eq:131} about $\bm s_0=\bm e$ can be simplified by considering $\left\|\bm 0\right\|_1=\left\|\bm 0\right\|_2$ as:
\begin{equation}\begin{array}{*{20}{l}}
{\begin{array}{*{20}{l}}
{{{\left\| {{s_k}} \right\|}^{\sum\limits_{vk = 1}^{n - k - 1 - j} {{b_{vk}}\left[ {{m_k} - 2{h_p} + \sum\limits_{vk = 1}^{n - k - 1 - j} {{b_{vk}}} (1 + {b_1}(1 - {m_k}))} \right]} }}{s_k}^{(j - l)}}\\
{ \cdot {{\left\| {{s_k}} \right\|}^{ - {m_k} - 2h}}{{\left\| {{s_k}} \right\|}^{\sum\limits_{vk = 1}^l {{b_{vk}}(1 + {b_1}(1 - {m_k}))} }}}
\end{array}}\\
{\begin{array}{*{20}{l}}
{\mathop {{\rm{  }} = }\limits^{e \to 0} {{\left\| {{s_k}} \right\|}^{{m_k}({h_p}(1 - {b_1}{h_p}) - 1 - h{b_1})}}}\\
{ \cdot {{\left. {{{\left\| {{s_k}} \right\|}^{(h_p^2 + h)({b_1} - 1) + 1}}{E_{n1}}} \right|}_{{h_v} = 1,2}}}
\end{array}}
\end{array}\label{eq:135}\end{equation}
with
\begin{equation}{\left\{ {\begin{array}{*{20}{l}}
{\sum\limits_{vk = 1}^{n -k- 1 - j} {{b_{vk}}} } = {h_{p}} \le {n-k-1-j}\\
~~~{\sum\limits_{vk = 1}^{l} {{b_{vk}}} } = h \le {l}\\
~~~~~{{h_{v}}- ({b_2} +  \cdots  + {b_{{l_v}})} = {b_1}} \le {l_v}
\end{array}} \right.}\label{eq:136}\end{equation}
where $j=l$ because the fianl exponential term of $\bm e^{(j-l)}$ will be smaller than $1$ if $j\ne l$ by considering the decomposition term $\sum\limits_{k = 0}^{n - 2} {(\frac{1}{{{T_c}}}{{\left\| {\frac{{\partial W({{\left\| {{\bm s_k}} \right\|}^{m_{k}}})}}{{\partial {\bm s_k}}}} \right\|}^{ - 2}}{{\left[ {\frac{{\partial W({{\left\| {{\bm s_k}} \right\|}^{m_{k}}})}}{{\partial {\bm s_k}}}} \right]}^T}} {)^{(n - k - 1)}}$ in \eqref{eq:40}.

Then, \eqref{eq:135} will not be infinite if the following inequality is satisfied:
\begin{equation}\begin{array}{l}
{m_k}( - {h_p}(1 - {b_1}{h_p}) + 1 + h{b_1}) \le (h_p^2 + h)({b_1} - 1) + 1\\
 \Rightarrow \left\{ {\begin{array}{*{20}{c}}
{{m_k} \le (\frac{{(h_p^2 + h) + 1}}{{({h_p}(2{h_p} - 1) + 1 + 2h)}})}\\
{{m_k} \le (\frac{{1}}{{({h_p}({h_p} - 1) + 1 + h)}})}
\end{array}} \right.\\
 \Rightarrow {m_k} \le (\frac{1}{{(1(1 - 1) + 1 + n -k- 1)}}) \le \frac{1}{n-k}
\end{array}\label{eq:137}\end{equation}
where $m_k$ will be minimum by considering the inequality relationship in \eqref{eq:136}.

Then the controller in condition of step 0 will not be singular if $0<m_k\le \frac{1}{n-k}$

Step $k+1,\cdots,n-2$:~~~~~~~~~~~\vdots

Step $n-1$: it is cleat that the singular term of controller about $s_{n-1}$ is the final term of controller in \eqref{eq:124}, and the singularity will not occur if $0<m_{n-1}<1$ by considering Lemma 2.

According to the these steps, a selection scheme of nonsingular parameters $m_k$ can be designed as $0<m_k\le \frac{1}{n-k}$ if Assumption 1 is satisfied.

The proof is completed.

According to Result 2 and Result 3, a final selection scheme of nonsingular parameters $m_k$ which can guarantee the predefined-time can be designed as \eqref{eq:112} and \eqref{eq:113}.

\section{Simulation and Analysis}

In order to show the effective of the proposed nonsingular sliding control method, two various-order control system will be shown in the following. In addition, the curve formed by singular parameters $m_k$ cannot be shown, because singularity is inevitable. Even if it is simulated by luck, it will not pass the  above-mentioned mathematical argument.

\textbf{Example 2:} Fig.4 shows a bad case: $m=0.9$ in second-order system and $m=0.5$ in third-order system. It is clear that the control input is infinite or very large in 1.143s in Fig.4 (1), and the control input is also very large in 0.805s in Fig.4 (2).

\begin{figure}
\begin{center}
\includegraphics[height=5cm]{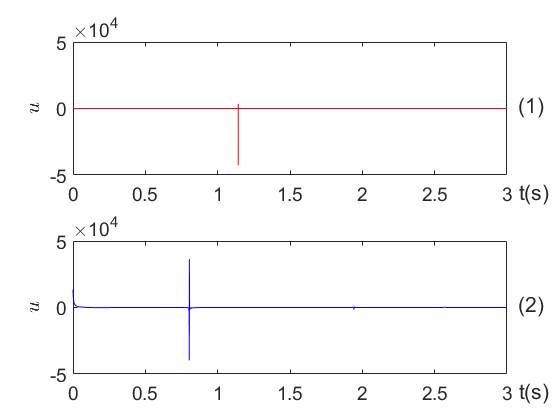}    
\caption{Some examples: (1)The curve of $y=u$ in the $W_1(x)$.$T_c=1$, $m=0.9$~~~(2)The curve of $y=u$ in the $W_1(x)$.}

{$T_c=1$, $m=0.5$. $cos(t)$ is the tracking signal}
\label{fig4}                                 
\end{center}                                 
\end{figure}

\subsection{A Second-order MIMO System}
\begin{equation}\left\{ {\begin{array}{*{20}{l}}
{\dot {\bm x_1} = {\bm x_2}}\\
{\dot {\bm x_2} =  - {\bm x_2} + \bm u}\\
\bm y=\bm x_1
\end{array}} \right.\label{eq:201}\end{equation}
where $\bm x_{1d}=\left[ {\begin{array}{*{20}{c}}
{\cos (t)}&1&{ - 5}&5
\end{array}} \right]^T$, the tracking error is $\bm e=\bm x_1-\bm x_{1d}$

A predefined-time function can be chosen as: $W_1(\bm x)=sin(arctan(\left\|\bm x\right\|^m))$. Then, a nonsingular predefined-time terminal sliding mode surface is designed by Corollary 2 as the following:
\begin{equation}\bm s = \bm e + \frac{1}{{{T_c}}}{\left[ {\frac{{\partial W_1({{\left\| \bm e \right\|}^m})}}{{\partial \bm e}}} \right]^{ - 1}}\label{eq:202}\end{equation}
with
\begin{equation}m = \left\{ {\begin{array}{*{20}{c}}
{1 + (m - 1)sign(sat(\frac{{\left\|\bm  e \right\|}}{\varepsilon })){\rm{~,~~~ }}\left\| \bm s \right\| \ne 0}\\
{{\rm{       ~~~~~~~~~~~~~~~~~           }}m{\rm{       ~~~~~~~~~~~~~~~~~~~,~~~         }}\left\| \bm s \right\| = 0}
\end{array}} \right..\label{eq:203}\end{equation}
where $T_c=1$, $\varepsilon=0.0001$, $u_{max}=15$.

Then, according to \eqref{eq:40}, the control input $\bm u$ can be designed as:
\begin{equation}\begin{array}{l}
\bm u =    {\bm x_2} - \sum\limits_{k = 0}^{n - 2} {\frac{1}{{m{T_c}}}({{\left[ {\frac{{\partial W({{\left\| {{{\bm s}_k}} \right\|}^m})}}{{\partial {{\left\| {{{\bm s}_k}} \right\|}^m}}}} \right]}^{ - 1}}} \frac{{{{\bm s}_k}}}{{{{\left\| {{{\bm s}_k}} \right\|}^m}}}{)^{(n - k - 1)}}\\
~~~~~~ { - \frac{1}{{m{T_c}}}{{\left[ {\frac{{\partial W({{\left\| {{{\bm s}_{n - 1}}} \right\|}^m})}}{{\partial {{\left\| {{{\bm s}_{n - 1}}} \right\|}^m}}}} \right]}^{ - 1}}\frac{{{{\bm s}_{n - 1}}}}{{{{\left\| {{{\bm s}_{n - 1}}} \right\|}^m}}}}
\end{array}\label{eq:204}\end{equation}
where $n=2$.

\begin{figure}
\begin{center}
\includegraphics[height=5cm]{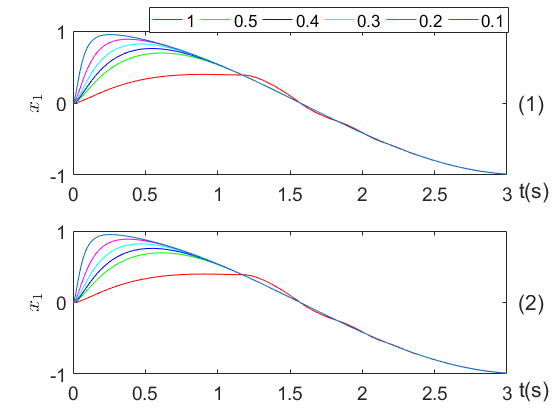}    
\caption{{(1)the method in \cite{ni2019predefined}}
{(2)the proposed method}}
{~~~~~~~~~~~~~The curve of $y=x_1$ in the $W_1(x)$.}
{$T_c=1$, $m=1,0.5,\cdots 0.1$} 
\label{fig5}                                 
\end{center}                                 
\end{figure}

\begin{figure}
\begin{center}
\includegraphics[height=5cm]{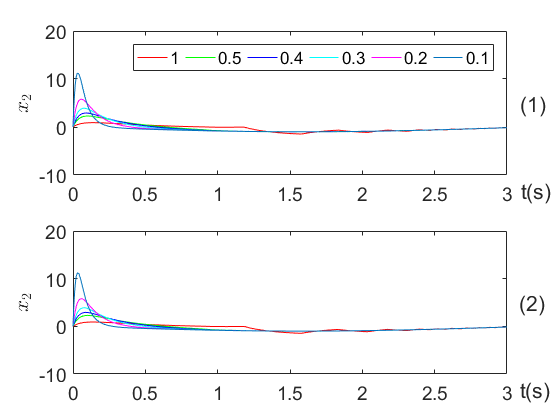}    
\caption{{(1)the method in \cite{ni2019predefined}}
{(2)the proposed method}}
{~~~~~~~~~~~~~The curve of $y=x_2$ in the $W_1(x)$.}
{$T_c=1$, $m=1,0.5,\cdots 0.1$}   
\label{fig6}                                 
\end{center}                                 
\end{figure}

\begin{figure}
\begin{center}
\includegraphics[height=5cm]{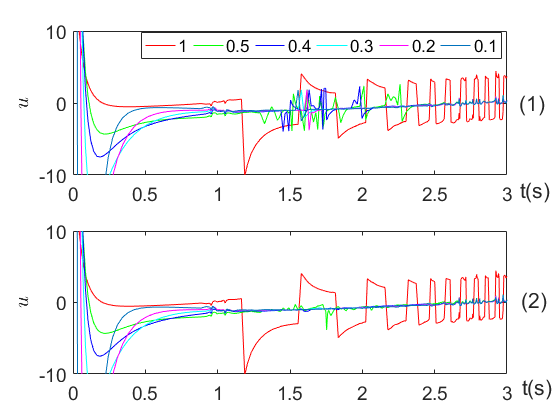}    
\caption{{(1)the method in \cite{ni2019predefined} (2)the proposed method}}
{~~~~~~~~~~~~~The curve of $y=u$ in the $W_1(x)$.}
{$T_c=1$, $m=1,0.5,\cdots 0.1$}
\label{fig7}                                 
\end{center}                                 
\end{figure}

\begin{figure}
\begin{center}
\includegraphics[height=5cm]{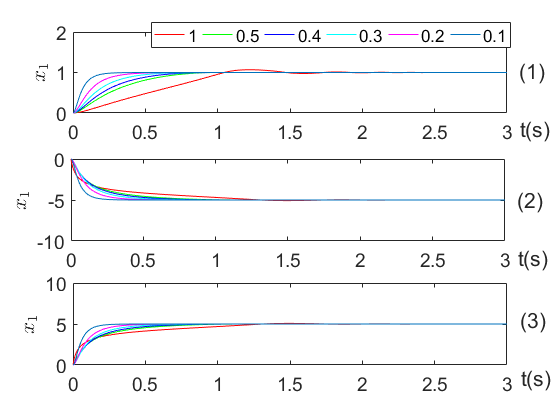}    
\caption{{(1)$x_{1d}=1$ ~~~~~~~(2)$x_{1d}=-5$ ~~~~~~~(3)$x_{1d}=5$}}
{~~~~~~~~~~~~~~~~~The curve of $y=x_1$ in the $W_1(x)$.}
{$T_c=1$, $m=1,0.5,\cdots 0.1$} 
\label{fig8}                                 
\end{center}                                 
\end{figure}

\begin{figure}
\begin{center}
\includegraphics[height=5cm]{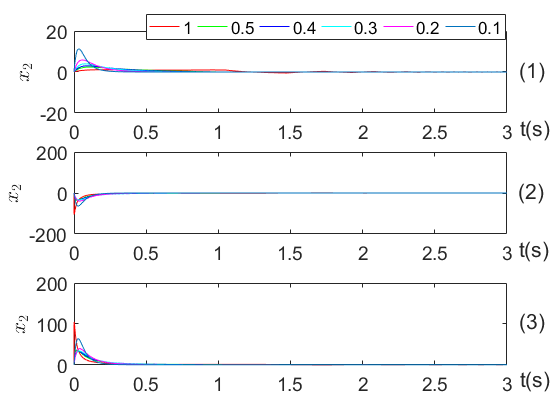}    
\caption{{(1)$x_{1d}=1$ ~~~~~~~(2)$x_{1d}=-5$ ~~~~~~~(3)$x_{1d}=5$}}
{~~~~~~~~~~~~~~~~~The curve of $y=x_1$ in the $W_1(x)$.}
{$T_c=1$, $m=1,0.5,\cdots 0.1$}  
\label{fig9}                                 
\end{center}                                 
\end{figure}

According to Introduction, the case of Remark 4 was not considered in the method in \cite{ni2019predefined}, so the unnecessary buffeting will occur in control input caused by inappropriate singular risk parameter $\varepsilon $. As shown in Figs.5 and 6, it is clear that both methods can make the tracking error converge to zero within predefined-time $2T_c=2s$ if $m$ is satisfied the nonsingular parameter range deduced by this paper. As shown in Fig 7, it is also clear that the chattering of the proposed control input is lower than the controller in \cite{ni2019predefined}. The tracking errors of second-order MIMO system can be shown in Figs. 8 and 9, and the tracking error $\bm e$ can not converge to zero within predefined-time if $m=1$.

\subsection{A Third-order MIMO System}
\begin{equation}\left\{ {\begin{array}{*{20}{l}}
{\dot {\bm x_1} = {\bm x_2}}\\
{\dot {\bm x_2} = {\bm x_3}}\\
{\dot {\bm x_3} =- {\bm x_2}  - {2\bm x_3} + \bm u}\\
\bm y=\bm x_1
\end{array}} \right.\label{eq:205}\end{equation}
where $\bm x_{1d}=\left[ {\begin{array}{*{20}{c}}
{\cos (t)}&1&{ - 5}&5
\end{array}} \right]^T$, the tracking error is $\bm e=\bm x_1-\bm x_{1d}$

A predefined-time function can be chosen as: $W_1(\bm x)=sin(arctan(\left\|\bm x\right\|^{m_k}))$. Then, a nonsingular predefined-time terminal sliding mode surface is designed by Corollary 2 as the following:
\begin{equation}\bm s_{k+1} = \bm s_k + \frac{1}{{{T_c}}}{\left[ {\frac{{\partial W_1({{\left\| \bm s_k \right\|}^{m_k}})}}{{\partial \bm s_k}}} \right]^{ - 1}}\label{eq:206}\end{equation}
with
\begin{equation}m_k = \left\{ {\begin{array}{*{20}{c}}
{1 + (m_k - 1)sign(sat(\frac{{\left\|\bm  s_k \right\|}}{\varepsilon })){\rm{~,~~~ }}\left\| \bm s_{k+1} \right\| \ne 0}\\
{{\rm{       ~~~~~~~~~~~~~~~~~~           }}m_k{\rm{       ~~~~~~~~~~~~~~~~~~~~,~~~         }}\left\| \bm s_{k+1} \right\| = 0}
\end{array}} \right..\label{eq:207}\end{equation}
where $k=0,1,2$, $\bm s_3=0$, $\bm s_0=\bm e$, $T_c=5$, $\varepsilon=0.001$, $u_{max}=50$. For convenience, $m_k=m=1,0.3,0.2,0.1$.

Then, according to \eqref{eq:40}, the control input $\bm u$ can be designed as:
\begin{equation}\begin{array}{l}
\bm u =    {\bm x_2}  - {2\bm x_3} - \sum\limits_{k = 0}^{n - 2} {\frac{1}{{m{T_c}}}({{\left[ {\frac{{\partial W({{\left\| {{{\bm s}_k}} \right\|}^m})}}{{\partial {{\left\| {{{\bm s}_k}} \right\|}^m}}}} \right]}^{ - 1}}} \frac{{{{\bm s}_k}}}{{{{\left\| {{{\bm s}_k}} \right\|}^m}}}{)^{(n - k - 1)}}\\
~~~~~~ { - \frac{1}{{m{T_c}}}{{\left[ {\frac{{\partial W({{\left\| {{{\bm s}_{n - 1}}} \right\|}^m})}}{{\partial {{\left\| {{{\bm s}_{n - 1}}} \right\|}^m}}}} \right]}^{ - 1}}\frac{{{{\bm s}_{n - 1}}}}{{{{\left\| {{{\bm s}_{n - 1}}} \right\|}^m}}}}
\end{array}\label{eq:208}\end{equation}
where $n=3$.

\begin{figure}
\begin{center}
\includegraphics[height=5cm]{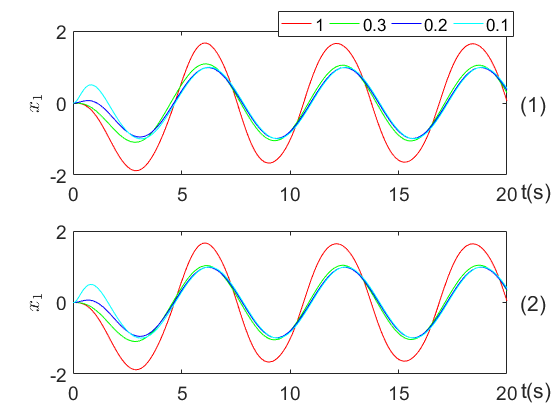}    
\caption{{(1)the method in \cite{ni2019predefined} (2)the proposed method}}
{~~~~~~~The curve of $y=x_1$ in $W_1(x)$}
{$T_c=5$, $m=1,0.3,\cdots 0.1$}  
\label{fig10}                                 
\end{center}                                 
\end{figure}

\begin{figure}
\begin{center}
\includegraphics[height=5cm]{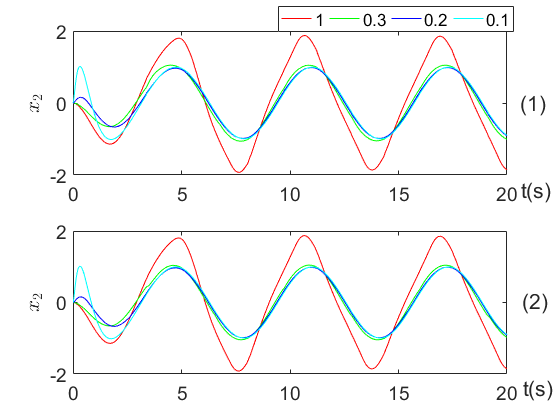}    
\caption{{(1)the method in \cite{ni2019predefined} (2)the proposed method}}
{~~The curve of $x_2$ in $W_1(x)$}
{$T_c=5$, $m=1,0.3,\cdots 0.1$}  
\label{fig11}                                 
\end{center}                                 
\end{figure}

\begin{figure}
\begin{center}
\includegraphics[height=5cm]{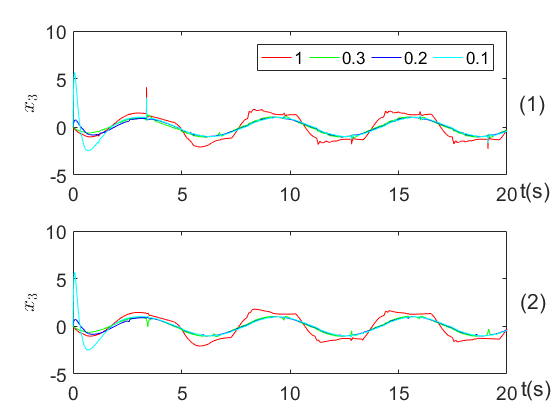}    
\caption{{(1)the method in \cite{ni2019predefined} (2)the proposed method}}
{~~The curve of $u$ in $W_1(x)$}
{$T_c=5$, $m=1,0.3,\cdots 0.1$}  
\label{fig12}                                 
\end{center}                                 
\end{figure}

\begin{figure}
\begin{center}
\includegraphics[height=5cm]{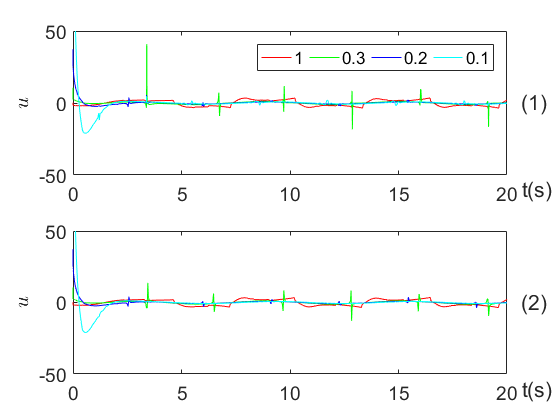}    
\caption{{(1)the method in \cite{ni2019predefined} (2)the proposed method}}
{~~~~~~~~The curve of $\bm y=\bm x_1$ in $W_1(x)$}
{$T_c=5$, $m=1,0.3,\cdots 0.1$} 
\label{fig13}                                 
\end{center}                                 
\end{figure}

\begin{figure}
\begin{center}
\includegraphics[height=5cm]{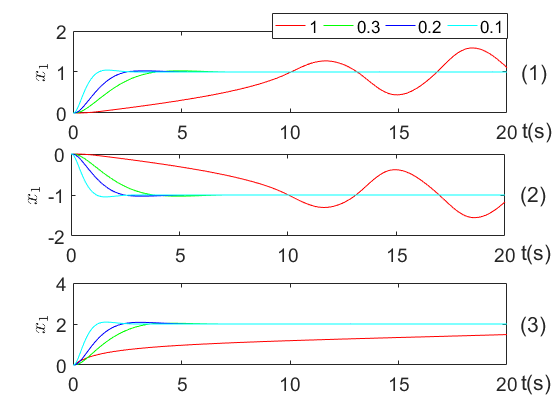}    
\caption{{(1)$x_{1d}=1$ ~~~~~~~(2)$x_{1d}=-1$ ~~~~~~~(3)$x_{1d}=2$}}
{~~~~The curve of $\bm x_2$ in $W_1(x)$}
{$T_c=5$, $m=1,0.3,\cdots 0.1$}  
\label{fig14}                                 
\end{center}                                 
\end{figure}

\begin{figure}
\begin{center}
\includegraphics[height=5cm]{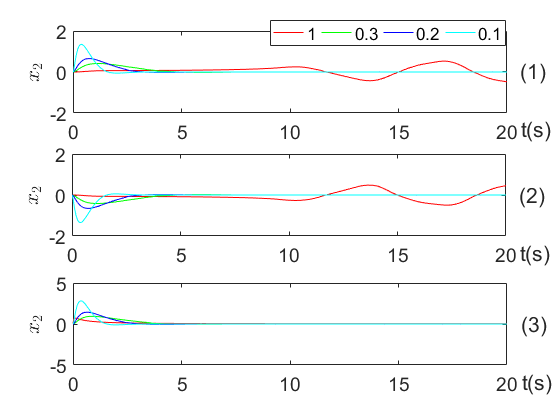}    
\caption{{(1)$x_{1d}=1$ ~~~~~~~(2)$x_{1d}=-1$ ~~~~~~~(3)$x_{1d}=2$}}
{~~~~The curve of $\bm x_2$ in $W_1(x)$}
{$T_c=5$, $m=1,0.3,\cdots 0.1$} 
\label{fig15}                                 
\end{center}                                 
\end{figure}

\begin{figure}
\begin{center}
\includegraphics[height=5cm]{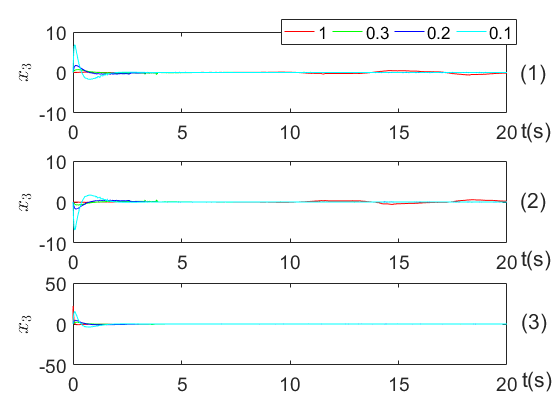}    
\caption{{(1)$x_{1d}=1$ ~~~~~~~(2)$x_{1d}=-1$ ~~~~~~~(3)$x_{1d}=2$}}
{~~~~The curve of $\bm x_2$ in $W_1(x)$}
{$T_c=5$, $m=1,0.3,\cdots 0.1$}  
\label{fig16}                                 
\end{center}                                 
\end{figure}

According to Introduction, the case of Remark 4 was not considered in the method in \cite{ni2019predefined}, so the unnecessary buffeting will also occur in control input caused by inappropriate singular risk parameter $\varepsilon $. As shown in Figs.10, 11 and 12, it is clear that both methods can make the tracking error converge to zero within predefined-time $2T_c=2s$ if $m$ is satisfied the nonsingular parameter range deduced by this paper. As shown in Fig 13, it is also clear that the chattering of the proposed control input is also lower than the controller in \cite{ni2019predefined}. The tracking errors of third-order MIMO system can be shown in Figs. 14, 15 and 16, and the tracking error $\bm e$ can not converge to zero within predefined-time if $m=1$.

\section{Conclusion}
Generally, the terminal sliding mode of fixed-time-class is often accompanied by the inevitable singularity problem and many researches have been done to overcome the singularity, but often at the larger cost of some controller performance. In this paper, we discuss a selection problem of nonsingular parameter in a predefined-time terminal sliding mode control strategy in as much detail as possible. Firstly, we discuss the parameter selection of traditional general predefined-time system expression equation and improve some conditions to make the system have a better stable performance. Secondly, we study the selection problem of nonsingular parameter for second-order system: a switching strategy is improved to overcome the singularity of reaching stage of sliding mode and chattering caused by switching sliding mode (Remark 4), and then, a selection range of nonsingular parameter is found to ensure the non-singularity of the sliding stage(s=0). Thirdly, we try to choose a recursive terminal sliding mode control method to deal with the predefined-time sliding mode control problem of high-order system, in this subsection, a selection method of nonsingular parameter is studied to try to achieve nonsingular predefined-time terminal sliding mode control for a n-order system. In the end, we use some simple system model to carry on the relevant numerical simulation, the final results, to some extent, confirm our inference.
         \bibliographystyle{autart}
\bibliography{mybibfile}                                  



\end{document}